\documentclass[preprint]{elsarticle}
\usepackage[utf8]{inputenc}
\usepackage{subcaption}
\usepackage{graphicx}
\usepackage{grffile}
\usepackage{epstopdf}
\usepackage{amsmath}
\usepackage{amssymb}
\usepackage{bm}
\usepackage{hyperref}
\usepackage{xcolor}
\usepackage[autostyle]{csquotes}
\usepackage{listings}
\usepackage{todonotes}
\usepackage{array}

\journal{Elsevier}

\newcommand{\sym}[1]{(\protect\includegraphics[height=6pt]{styles/#1.eps})}
\newcommand{\syms}[1]{\protect\includegraphics[height=6pt]{styles/#1.eps}}

\begin{document}

\begin{frontmatter}

\title{Thermodynamically consistent physics-informed neural networks for hyperbolic systems}
\author[1]{Ravi G. Patel}
\ead{rgpatel@sandia.gov}
\author[4]{Indu Manickam}
\ead{imanick@sandia.gov}
\author[1]{Nathaniel A. Trask$^*$}
\ead{natrask@sandia.gov}
\author[2]{Mitchell A. Wood}
\ead{mitwood@sandia.gov}
\author[3]{Myoungkyu Lee}
\ead{mnlee@sandia.gov}
\author[1]{Ignacio Tomas}
\ead{itomas@sandia.gov}
\author[1]{Eric C. Cyr}
\ead{eccyr@sandia.gov}

\cortext[cor1]{Corresponding author}
\address[1]{Sandia National Laboratories - Computational Mathematics Department}
\address[2]{Sandia National Laboratories - Computational Multiscale Department}
\address[3]{Sandia National Laboratories - 
Combustion Research Facility}
\address[4]{Sandia National Laboratories - 
Mission Algorithms Research \& Solutions}

\begin{abstract}

Physics-informed neural network architectures have emerged as a powerful tool for developing flexible PDE solvers which easily assimilate data, but face challenges related to the PDE discretization underpinning them. By instead adapting a least squares space-time control volume scheme, we circumvent issues particularly related to imposition of boundary conditions and conservation while reducing solution regularity requirements. Additionally, connections to classical finite volume methods allows application of biases toward entropy solutions and total variation diminishing properties. For inverse problems, we may impose further thermodynamic biases, allowing us to fit shock hydrodynamics models to molecular simulation of rarefied gases and metals. The resulting data-driven equations of state may be incorporated into traditional shock hydrodynamics codes.

\end{abstract}


\end{frontmatter}


\section{Introduction}


Recently, a number of works have evaluated the potential of deep neural networks (DNNs) to solve partial differential equations (PDEs). DNNs possess attractive properties: potentially exponential convergence, breaking of the curse-of-dimensionality, and an ability to handle data sampled from function spaces with limited regularity, such as shock and contact discontinuities \cite{weinan2018deep,he2018relu,daubechies2019nonlinear,yarotsky2017error,yarotsky2018optimal,opschoor2019deep,bach2017breaking,bengio2000taking,han2018solving}. Practically however, challenges regarding the training of DNNs often prevent the realization of convergent schemes for forward problems \cite{wang2020understanding,beck2019full,fokina2019growing,adcock2020gap}. For inverse problems however, a number of methods have emerged that train neural networks to simultaneously match target data while minimizing a PDE residual \cite{lagaris1998artificial,raissi2019physics,raissi2018deep}, which have found application across a wide range of applied mathematics problems \cite{sun2020surrogate,zhang2019quantifying,meng2020composite,mao2020physics,zhang2020learning}. While many variations of this idea exist in the literature, we develop in this work an extension to physics-informed neural networks (PINNs) applied to hyperbolic systems of conservation laws. For this method and related PINNs approaches, the low-dimensional neural network representation of a PDE solution allows a simple and efficient implementation of inverse problems in popular machine learning packages such as Tensorflow~\cite{tensorflow2015-whitepaper} and PyTorch~\cite{NEURIPS2019_9015}. 

We consider the application of these approaches to inverse problems in shock physics. Specifically we consider shock-hydrodynamics, a class of models where the Euler equations are used to model materials in high-energy regimes under which shear may be assumed negligible (i.e. substances subject to very high temperatures and rates of deformation), and may be extended to include elastoplastic material response in addition to electromagnetic physics \cite{robinson2013fundamental}. For these problems, an accurate and appropriate equation of state (EOS) is critical as the material of interest evolves over a large area of phase-space including phase transitions. The process of developing an EOS is typically handled by a labor-intensive assimilation of heterogeneous data compiled from analytical, experimental, and synthetic sources \cite{carpenter2015automated}. In this work, we develop improvements upon traditional PINNs to support the extraction of data-driven EOS. For this application, incorporation of a number of physical principles will prove critical to obtaining an EOS which not only fits training data well, but also provides stable solutions when incorporated into a traditional continuum shock solver after training.

The hybrid physics/data loss used in PINNs admits interpretation as a weighted least-squares collocation scheme. As such, PINNs inherit a number of disadvantages of such methods: a need to appropriately weight PDE residual against initial and boundary conditions, a restrictive regularity requirement that solutions be continuous, and a lack of natural means to enforce conservation structure. As a result, applications of PINNs to solve hyperbolic forward problems requires introduction of a number of penalties \cite{jagtap2020conservative,mao2020physics}. The premise of the current work is that an alternative space-time least-squares control volume discretization (cvPINNs) naturally resolves these issues, substantially removing hyper-parameters while providing higher quality solutions. At the same time, exposing connections to traditional finite volume methods (FVM) allows introduction of novel techniques to ensure the solution is asymptotically consistent with the zero-viscosity limit and satisfaction of physically relevant entropy inequalities \cite{lax1973hyperbolic,menikoff1989riemann}. In order to do so, we use the framework of complete-EOS \cite{Meni1989}, which allows us to write down an entropy-flux pair and enforce entropy conditions in the forward-model. For the inverse model, we parameterize the specific entropy describing the entire thermodynamic behaviour of our substance/material of interest. In this setting, it is important to enforce a number of inequality constraints \cite{Viscous2014} into the specific entropy in order to preserve, among other things, hyperbolicity in the Euler equations.


A number of works have considered the use of DNNs in the context of shock problems \cite{jagtap2020conservative,mao2020physics,xiong2020roenets,tokareva2019machine,magiera2020constraint}, and several works have considered using alternative discretizations in the context of PINNs-like methods, for example Ritz-Galerkin discretizations \cite{weinan2018deep}, Petrov-Galerkin methods \cite{kharazmi2020hp}, and mortar methods \cite{jagtap2020extended}. To our  knowledge, this work marks the first attempt to assimilate traditional finite volume methodology to obtain a thermodynamically consistent treatment of inverse problems in shock physics.

We organize the paper by providing an interpretation of classical conservation laws as an integral balance law in space-time, together with an inequality constraint on entropy. We then present both classical PINNs and cvPINNs in the context of forward modeling before moving toward inverse problems. We provide pedagogical examples for a variety of canonical hyperbolic systems, such as Burgers, Euler, and Buckley-Leverett equations. We consider the extraction of equations of state from realistic high-fidelity noisy synthetic sources. First, direct simulation Monte Carlo (DSMC) simulation of Argon gas in the continuum regime allows a comparison of fitting a general purpose DNN EOS vs. a perfect gas law. We conclude by extracting an EOS from molecular dynamics (MD) simulations of shocks propagating through copper bars thereby demonstrating the framework's ability to handle non-fluid materials whose EOS is \textit{a priori} unknown.

\section{Space-time Integral Form PDE Formulation}

Given a space-time domain $\Omega \subset \mathbb{R}^d \times [0,T]$ and a conserved vector quantity $\bm{u} \in \mathbb{R}^P$, we consider a class of conservation laws of the form
\begin{align}\label{eq:pde}
    \begin{split}
    \partial_t \bm{u} + \nabla \cdot \bm{F}(\bm{u}) = 0 &\qquad x,t \in \Omega, \text{ for all } i\\
    \bm{u} = \bm{u}_0 &\qquad t = 0\\
    \bm{F}(\bm{u})\cdot \hat{n} = g &\qquad x \in \Gamma_{-}
    \end{split}
\end{align}
where bold denotes a vector field, $\bm{F}\in \mathbb{R}^{d \times P}$ a flux, and $\Gamma_{-}$ a problem specific subset of $\partial \Omega$ with positive measure. We refer to components of $\bm{u}$ as $u^p$ where the superscript $p\in \left\{1,...,P\right\}$ is the $p^{th}$ vector component. $\Gamma_{-}$ is generally associated with the ``upwind'' portion of the boundary. 
Equation~\eqref{eq:pde}, interpreted in a weak sense, potentially leads to non-physical multi-valued solutions \cite{Darfermos2000}. Following \cite{Bianchi2005}, we are only interested in solutions understood as the \textit{zero-viscosity limit} given by $\bm{u} = \lim_{\epsilon \rightarrow 0^+}\bm{u}_\epsilon$ where
\begin{align}\label{eq:limpde}
        \partial_t \bm{u}_\epsilon + \nabla \cdot \bm{F}(\bm{u}_\epsilon) = \epsilon \nabla^2 \bm{u}_\epsilon ,
\end{align}
also called viscosity-solutions. Let $\eta:\mathbb{R}^N \rightarrow \mathbb{R}$ be a convex functional and $\bm{q}(\bm{u}):\mathbb{R}^N \rightarrow \mathbb{R}^{d}$ a vector-valued function satisfying the identity $\nabla_{\bm{u}}\eta(\bm{u})^{\top} \nabla_{\bm{u}} [\bm{F}(\bm{u})^{\top}] = \nabla_{\bm{u}}\bm{q}(\bm{u})$.  Then $(\eta(\bm{u}),\bm{q}(\bm{u}))$ is called an \textit{entropy-flux pair}. Viscosity solutions satisfy the so-called entropy inequality (see for instance \cite[p. 21]{GodRav1996})
\begin{equation}\label{eq:pdie}
    \partial_t \eta + \nabla \cdot \bm{q} \leq  0 ,
\end{equation}
with equality holding only in smooth regions of the solution. The convex functional $\eta$ is known as a \textsl{mathematical entropy}. We highlight that this is a different object from $s$, the \textsl{specific entropy}, to be introduced later in this manuscript.

We define an ``extended-flux'' $\hat{\bm{F}} := \left< \bm{u}^\intercal, \bm{F}\right> \in \mathbb{R}^{d+1 \times P}$ which has the additional column $\bm{u}^\intercal$. Similarly, we may define an extended entropy-flux to obtain the space-time entropy flux $\hat{\bm{q}} = \left< \eta, \bm{q}\right> \in \mathbb{R}^{d+1}$. Thus the conservation law and entropy inequality may be written in terms of the generalized divergences
\begin{equation}\label{spacetimeClaw} 
    \begin{aligned}
&div(\hat{\bm{F}}) = 0, \\
&div(\hat{\bm{q}}) \leq 0 ,
\end{aligned}
\end{equation}
where $div = \left< \partial_t, \partial_{x_1},...,\partial_{x_d}\right>$ to explicitly distinguish from the space-only $\nabla \cdot$ operator. Application of the Gauss divergence theorem yields the following integral conservation law form, which holds for any compact $\omega \subset \Omega$ with piecewise smooth boundary $\partial \omega$,
\begin{equation}\label{spacetimeClawint}
    \begin{aligned}
&\int_{\partial \omega} \hat{\bm{F}} \cdot d\bm{A} = 0, \\
&\int_{\partial \omega} \hat{\bm{q}} \cdot d\bm{A} \leq 0. 
\end{aligned}
\end{equation}

\section{Physics-informed neural networks}

To emphasize the advances developed here, this section recalls a classical view of PINNs as a point collocation least squares method and develops a new control volume PINNs (cvPINNs) formulation that has notable advantages for hyperbolic conservation laws. 

\subsection{Classical PINNs: Point collocation least squares}

Denote by $|| f ||_{\ell_2(\mathcal{D})}$ the root-mean-square norm of a collection of scattered point data $\mathcal{D} = \left\{\bm{x}_i,\bm{f}(\bm{x}_i)\right\}_{i=1,...,N_{data}}$. The classical PINNs approach to solve Eqn. \ref{eq:pde} introduces pointsets on the interior ($\mathcal{D}_{int} \subset \Omega$), boundaries ($\mathcal{D}_{BC}\subset \Gamma_{-}$), and initial time points ($\mathcal{D}_{IC}\subset \Omega \cap \left\{t=0\right\}$) to define the least squares residual
\begin{equation}\label{eq:oldpinn}
\begin{aligned}
\mathcal{L}_{\text{PINN}} = &\\
    ||\partial_t \bm{u} &+ \nabla\cdot \bm{F}(\bm{u})||_{\ell_2(\mathcal{D}_{int})}^2
   + \epsilon_0 ||\bm{u}-\bm{u}_0||_{\ell_2(\mathcal{D}_{IC})}^2
   + \epsilon_\Gamma ||\bm{F}(\bm{u})\cdot\hat{n}-g\||_{\ell_2(\mathcal{D}_{BC})}^2,
\end{aligned}
\end{equation}
where $\epsilon_0$ and $\epsilon_\Gamma$ are penalty hyperparemeters requiring calibration. The solution $\bm{u}$ is assumed to be a neural network whose solution is typically obtained by minimizing Eqn. \ref{eq:oldpinn} with first-order optimizers available in popular machine learning libraries (e.g., \cite{NEURIPS2019_9015,abadi2016tensorflow}). 

In the context of traditional PDE discretization, this class of weighted residual methods has been used extensively with choices of approximation other than DNNs \cite{moritz1978least,rummel1979least,ling2005least,hu2007weighted,zhang2001least,cheng2010collocation}, and possesses two fundamental challenges. First, the use of a point collocation functional mandates working in continuous function spaces, posing challenges for reduced regularity problems such as those occurring in shocks or $H(div)/H(curl)$ problems. Secondly, the proper weighting of penalty parameters is required to obtain coercivity of $\mathcal{L}_{\text{PINN}}$ over an appropriate energy norm. In least squares finite element contexts, tools such as the Agmon-Douglis-Nirenberg (ADM) theory \cite{BochevGLSbook} provides weights involving relative measures of cell volumes and boundary faces in the context of Elliptic theory. Neural networks however possess no explicit relationship to a static mesh; in fact DNNs may be identified with a piecewise linear finite element space which evolves during training \cite{he2018relu}, and therefore the requisite mesh information is not available a priori. These two challenges, stemming from the choice of weighted residual, lead to a number of pathologies in training PINNs that are an active area of current research \cite{wang2020understanding}. 


\subsection{Control Volume PINNs} \label{sec:cvpinns}
In control volume PINNs (cvPINNs) we take the space-time integral form, Eqn.~\ref{spacetimeClawint}, as the basis of generating a residual. We approximate $\bm{u} = \mathcal{NN}(t,x;\xi)$, where $\mathcal{NN}$ is a deep neural network taking $t$ and $\bm{x}$ as input and $\xi$ as parameters. Partitioning the domain $\Omega$ into $N_c$ disjoint space-time cells denoted $c$ and applying Eqn.~\ref{spacetimeClawint} provides the loss
\begin{equation}\label{eq:cvPINNloss}
\mathcal{L}_{cvPINN} = \sum_{c=1}^{N_c} \left| \int_{f \in \partial c} \widetilde{\bm{F}} \cdot d\bm{A} \right|^2.
\end{equation}
To naturally impose boundary conditions, the fluxes are evaluated conditionally as
\begin{equation}
    \widetilde{\bm{F}} = 
\begin{cases}
    \hat{\bm{F}}(\mathcal{NN}),& \text{if } \bm{x} \in \Omega\\
    g \hat{n},              & \text{if } \bm{x} \in \Gamma_{-}
\end{cases}
\end{equation}
that is, for cells touching an upwind boundary facet, the flux boundary condition is used instead. In this manner, the boundary and initial conditions are incorporated directly, eliminating the parameters $\epsilon_0$ and $\epsilon_\Gamma$. Conservation is naturally imposed via the control volume formulation, eliminating the need for conservation penalties (see e.g. \cite{Jagtap2020}).

In contrast to traditional FV methods, our approximation is not piecewise polynomial within cells with reconstructed piecewise polynomial fluxes at facets. Instead, $\mathcal{NN}|_f$ (the restriction of $\mathcal{NN}$ to $f$) inherits the regularity of the network activation functions; e.g. for a ReLU activation $\mathcal{NN}|_f$ is nonpolynomial and only $C_0$ continuous. Consequently, the quadrature in Eqn. \ref{eq:cvPINNloss} must be performed approximately. For this work we use either composite trapezoidal or composite midpoint quadrature, and a quadrature refinement study has been performed for all presented results to ensure that sufficiently many intervals have been used to reduce quadrature error below the optimization error in training.



We note that characterization of the solution of a hyperbolic system as the $L_2$-minimizer of the residual (as proposed in \eqref{eq:cvPINNloss}) might not necessarily retrieve the physically valid viscosity-solution. In other words, if viscosity-solutions of hyperbolic systems are meant to satisfy a Dirichlet principle (i.e. solutions can be characterized as minimizers of some residual) the right choice of norm for such residual is not known. This well-known fact has been studied by different authors. For instance in \cite{Bochev2016} the idea of adaptively weighted $L_2$-norms are advanced. For the case of scalar conservation laws and Hamilton-Jacobi equations $L_1$-minimization of the PDE-residual can be proven to retrieve the unique viscosity solution while the $L_2$-norm in general does not (see e.g. \cite{L1Guermond,Guermond2008}). For the case of systems, the right choice of norm is a largely open question.

Here we adopt a more pragmatic approach, choosing to work with the $L_2$-norm to allow use of generic optimization schemes available in popular machine learning packages. In order to mitigate potential deficiencies of the $L_2$-norm we introduce tools to both penalize the violations of the entropy inequality and total variation of our cvPINNs solution. These tools are critical to handling the discontinuities from shocks and contacts that lead to non-physical oscillations which may violate physical constraints. These constraints will be enforced by introducing penalties that allow a simple implementation. We stress however that unlike initial and boundary conditions, their weightings are not tied to a proper weighting by discretization lengthscales; later results will show that taking a unit weight for all results works well. The approaches presented here are not meant to represent the state-of-the-art of forward-solution methods for hyperbolic system of conservation laws; they are meant to indicate the types of FV tools which may be easily incorporated into physics-informed ML for shock problems.

For ease of presentation, we restrict the remainder of the section to a 1D spatial domain with a Cartesian mesh, although generalizations to polyhedral meshes are possible. We identify a space-time cell $c$ with centroid $(x_i,t_n)$ as $c_{i,n}$, and denote the north, east, south, and west facets as $f_{i,n+\frac12}$, $f_{i+\frac12,n}$, $f_{i,n-\frac12}$ and $f_{i-\frac12,n}$, respectively.




\paragraph{Artificial viscosity penalization} A traditional approach to obtaining entropy solutions is to introduce a mesh size-dependent viscosity which vanishes in the continuum limit, effectively discretizing Eqn. \ref{eq:limpde} by replacing $\epsilon$ with the characteristic mesh size $\Delta x$. We specifically consider the classical artificial viscosity from \cite{Reisner2013},
\begin{equation} \label{eq:avisc}
    \partial_t \bm{u} + \partial_x \mathbf{F}(\bm{u}) = \beta (\Delta x)^2 \partial_x \left( |\partial_x v|\, \partial_x \bm{u} \right)
\end{equation}
where $\Delta x$ is the cell width, $\beta>0$ is a small parameter, and $v$ denotes a problem-specific velocity field. This equation can be written in terms of the generalized space-time flux 
\begin{equation}\label{spacetimeClawVisc}
\hat{\mathbf{F}}_{AV} = \left<\bm{u}^\intercal,\bm{F}-\epsilon_{AV} (\Delta x)^2 \left( |\partial_x v| \partial_x \bm{u} \right)\right>
\end{equation}
 as above. In the remainder of the paper, taking $\epsilon_{AV}>0$ corresponds to adding in this additional flux.


While popular in classical FV methods, when working with DNN solutions the loss may only be incorporated to within optimization error. As a result, for sufficiently small $\Delta x$, there may be insufficient precision for the artificial viscosity to impact the optimizer. We thus introduce this as a means of comparison and motivation for the following alternative means of enforcing a vanishing-viscosity mechanism.


\paragraph{Entropy inequality penalization} Rather than working with Eqn.~\ref{eq:limpde}, we may instead seek to impose the entropy inequality constraint in Eqn.~\ref{spacetimeClawint} directly. While a formal treatment of inequality-constraints requires more sophisticated optimizers, we opt to penalize deviations from the constraint by implementing the following loss,
\begin{align}
    \label{eq:eres}
    &\mathcal{L}_{ent} = \sum_{c=1}^{N_c} \left(  \max\left(0,\int_{\partial c}\hat {\bm{q}}  \cdot d \bm{A} \right) \right)^2.
\end{align}

%


\paragraph{Total Variation penalization}
Treatment of discontinuous solutions in a variational $L_2$-framework inevitably leads to Gibbs-like phenomena. We highlight that, in the space-time setting of the cvPINNs solution, sub-optimal training/approximation of the DNN may account for part of these artificial oscillations. 


The total variation diminishing (TVD) property is given by
\begin{align}
\begin{split}\label{TVDbound}
    TV(\bm{u}_{n+1}) &\leq TV(\bm{u}_n), \\
    \text{where } TV(\bm{u}_n) &= sup \left\{\int_\Omega \bm{u}(\bm{x},t_n)\, div \phi\, d\bm{x}:\phi \in C^1_c, ||\phi||_{L^\infty} < 1\right\},
\end{split}
\end{align} 
where $C^1_c$ is the set of continuously differentiable, compactly supported vector functions on $\Omega$. Scalar conservation laws in one space dimension satisfy a total-variation bound, therefore the property \eqref{TVDbound} is a highly desirable feature in numerical solutions of scalar conservation laws (cf. \cite{GodRav1996}). However, exact solutions of neither scalar conservation laws in two-space dimensions nor general hyperbolic systems are guaranteed to satisfy a $TV$-bound (cf. \cite{Darfermos2000,Rauch1986}). In spite of this contradiction, introduction of a mild $TV$ penalization/regularization is a very popular numerical device used to mitigate spurious behaviour of the numerical solution (see e.g. \cite{Toro2000}). In the cvPINNs framework we use a similar approach.

%
In the classical FVM, $TV$ is often approximated on a 1D Cartesian mesh via
\begin{equation}
\begin{aligned}
    TV_h(\bm{u}_n) = \sum_i |\bm{u}_{i+1,n} - \bm{u}_{i,n}|,
\end{aligned}
\end{equation} 
where $\bm{u}_{i,n}$ are the states at the centers of cell, $(i,n)$. In one dimension this may be incorporated into the loss ,
\begin{equation}
    \mathcal{L}_{TVD} = \sum_{n,p} \max\left(0,  TV_h(u^p_{n+1}) - TV_h(u^p_n)  \right)^2.
\end{equation}

\paragraph{Physics-informed cvPINN}

We finally regularize the cvPINN loss with the artificial viscosity, entropy inequality, and TVD penalties to obtain the following loss functional governing forward simulation problems

\begin{equation} \label{eq:pderesTVD}
    \mathcal{L}_{fwd} = \mathcal{L}_{cvPINN} + \epsilon_{ent} \mathcal{L}_{ent} + \epsilon_{TVD} \mathcal{L}_{TVD}
\end{equation}
where $\epsilon_{ent}$ and $\epsilon_{TVD}$ are positive penalty hyperparameters which we will demonstrate how to set. It is understood that $\mathcal{L}_{cvPINN}$ is modified to accommodate the artificial viscosity flux in Eqn. \ref{spacetimeClawVisc} when we take $\epsilon_{AV}>0$.


\section{Least squares control volume scheme for inverse problems} \label{sec:inv}

This section generalizes the loss in Eqn. \ref{eq:pderesTVD} to inverse problems. We consider a class of conservation laws taking the form Eqn. \ref{eq:pde}, but where the exact functional form of the flux $\bm{F}$ is partially unknown. We assume that it may be parameterized as $\bm{F}_\sigma$, for parameter $\sigma$. This encompasses estimation of unknown material properties such as viscosity, as well as the more challenging setting where a model term, such as a multiscale closure or equation of state, are unknown. In this setting the parameter $\sigma$ could correspond to either a selection of candidate models from a dictionary \cite{kaiser2018sparse}, or a neural network parameterization of missing physics \cite{lu2019deeponet}.

We assume scattered point observations, $ \bm{u}_{data} $, at a pointset, $ \mathcal{D} \in \Omega$, which may correspond to partial observations of components of the state variable. We define the inverse problem as
\begin{equation}\label{inverse}
\underset{\xi,\sigma}{\text{min}}\,\,\mathcal{L}_I,\qquad \mathcal{L}_I = \mathcal{L}_{fwd} + \epsilon_{data} ||\mathcal{NN}-\bm{u}_{data}||^2_{\ell_2(\mathcal{D})}.
\end{equation}
Note that the minimization is performed simultaneously over neural network for the states (the $\xi$ parameters) and the parameters for the flux $\sigma$.

Shock hydrodynamics provides a motivating example where the specific form of the EOS maybe \emph{a priori} unknown. We will consider both parameter estimation of a known EOS model form (e.g. $\gamma$ in the ideal gas law), and a black-box DNN model appropriate for materials in extreme energy environments undergoing phase transitions. We specialize in the remainder to this EOS estimation problem. Consider as conservation law the 1D Euler equations,
\begin{equation} \label{eq:euler}
    \partial_t \begin{bmatrix} \rho \\ M \\ E \end{bmatrix} + \partial_x \begin{bmatrix} M \\ M^2/\rho + p \\ (E+p)M/\rho \end{bmatrix} = 0,
\end{equation}
where $\rho$ is the density, $M$ is the momentum, $E = \rho(e+{\frac{1}{2}u^2})$ is the total energy density, $p$ is the pressure, $u=M/\rho$ is the velocity, and $e$ is the specific internal energy. Expression \eqref{eq:euler} requires an additional EOS to obtain closure that characterizes the microscopic state of the system. 

While formula for the pressure $p=p(\rho,e)$ would suffice for closing \eqref{eq:euler}, we choose to parameterize a so-called \textit{complete-EOS}. A complete-EOS consists of an equation for the specific entropy $s := s(\rho,e)$, from which the temperature and pressure can be computed from the Gibbs relations:
\begin{align*}
T^{-1} = \frac{\partial s}{\partial e} 
\ \ \text{and} \ \
p = -\rho^2 \tfrac{\partial s}{\partial \rho}/ \tfrac{\partial s}{\partial e} \, .
\end{align*}
The framework of complete-EOS is particularly well described in the classical paper \cite{Meni1989} and has been used recently in \cite{Viscous2014} to establish the minimum principle of the specific entropy with largest generality possible. For our problem of interest, the description of the specific entropy has to be parametrized as $s := s(\rho,e;\sigma)$ where $\sigma$ is the parametrization. The compete-EOS description allows us to have well-defined entropy pairs by choosing \cite{Harten1998}
\begin{equation} \label{eq:eulerEnt}
    \eta = -\rho s; \quad \bm{q} = - M s.
\end{equation}
This choice of parameterization exposes the entropy pair required described in \S\ref{sec:cvpinns} .


However, not any function $s(\rho,e;\sigma)$ will yield physically admissible behaviour. The following inequality constraints are necessary to guarantee physically admissible \cite{Viscous2014} solutions:
\begin{align}
\label{eq:EOSconst}
\partial_e s > 0 \, , \
\partial_e^2 s \leq 0 \ \text{and} \
\partial_\rho (\rho^2 \partial_\rho s) < 0 .
\end{align}
Violations of the conditions in \eqref{eq:EOSconst} can result in a loss of hyperbolicity in the Euler equations and violation of the miniumum principle of the specific entropy.


We consider two parameterizations of the EOS. The first is a perfect gas, 
\begin{equation} \label{eq:eosPgas}
    s = \log(e^{1/(\gamma-1)}/\rho),
\end{equation}
parameterized by the ratio of specific heats, $\gamma$. To handle potentially negative $\rho$ or $e$ during training, we stabilize via 
\begin{equation}\label{eq:eosIdeal}
    s = \log\left(\frac{\max(\epsilon,e)^{1/(\gamma-1)}}{\max(\epsilon,\rho)}\right),
\end{equation}
where $\epsilon$ is a small number. It is easy to see that this model satisfies Eqns.~\ref{eq:EOSconst} for $\gamma>0$, and therefore will always provide a well-posed model.

As a second parametrization we take $s = \mathcal{NN}(\rho,e;\sigma)$. Of particular concern in situations with scarce data, this ``black-box'' parametrization provides no mechanism to enforce Eqns.~\ref{eq:EOSconst}. As before, we enforce these inequality constraints via penalty added to the loss in Eq.~\ref{inverse}
\begin{align}
    \label{eq:resIR}
    &\mathcal{L}_{inv} = \mathcal{L}_{I} + \epsilon_{R} \mathcal{L}_{R}
\end{align}
where 
\begin{align}
    \label{eq:resR}
    \begin{split}
    \mathcal{L}_{R} &= ||\max(0,-\partial_e s)||_{\ell_2(\mathcal{D}_{EOS})}^2 \\ 
    &+ ||\max(0,\partial_e^2 s)||_{\ell_2(\mathcal{D}_{EOS})}^2 \\ 
    &+ ||\max(0,\partial_\rho (\rho^2 \partial_\rho s))||_{\ell_2(\mathcal{D}_{EOS})}^2. 
    \end{split}
\end{align}
Here $\mathcal{D}_{EOS} = \{\rho_i,e_i\}^{i=1\hdots N_{EOS}}$ is a set of $N_{EOS}$ collocation points where the penalization is applied. In practice we identify a region of interest in parameter space and uniformly sample points. Together, the following three EOS parametrization (perfect gas parametrization, the neural network parametrization with entropy constraint penalization, and the unpenalized neural network) provide a sequence of increasingly ``black-box'' parametrization of the missing physics, and will serve as a test-bed for characterizing the role of physical inductive biases on generalizability and model stability in small data limits.

\section{Results: Forward solution of hyperbolic problems}

Before tackling the inverse problem discussed in Section~\ref{sec:inv}, we demonstrate that cvPINNs provides an appropriate discretization for nonlinear hyperbolic PDEs. 

\subsection{Entropy condition} \label{sec:ent}

\begin{figure}[htpb]
    \centering
    \includegraphics[width=\linewidth]{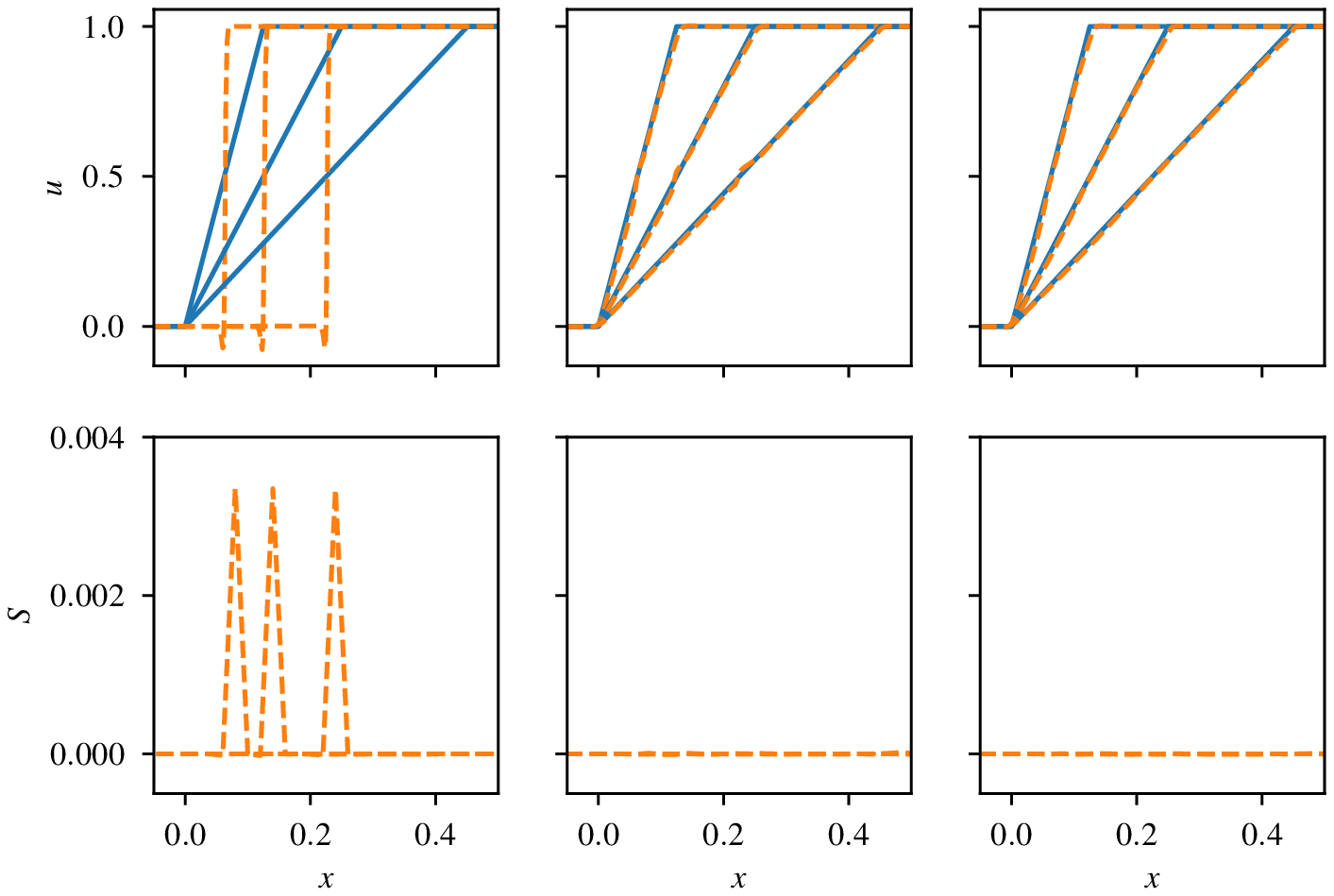}
    \caption{Effect of entropy regularization on cvPINNs solution for Burgers rarefaction Riemann problem with neural network $u$ initialized to rarefaction shocks. (\textit{top row}) cvPINNs solutions for TVD regularization weights, $\epsilon_{ent}=0$ (\textit{first column}), $\epsilon_{ent}=0.01$ (\textit{middle column}), and $\epsilon_{ent}=1$ (\textit{last column}). Analytical solution \sym{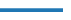} and cvPINNs solutions \sym{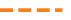} are shown for times, $t=0.125, 0.25, 0.45$. (\textit{bottom row}) Residual of the entropy inequality for the corresponding cvPINNs solutions. Without entropy penalization, training is stable with respect to the nonphysical rarefaction shock solution.}
    \label{fig:ent}
\end{figure}

The solution to Eqn.~\ref{eq:pde} is not unique without simultaneous satisfaction of Eqn.~\ref{eq:pdie}.  For instance, a Riemann problem setup for Burger's equations,  when  the right state exceeds the left (e.g. $u_L<u_R$), we might recover a nonphysical rarefaction shock solution if only the residual is minimized. 

To demonstrate this issue, and how entropy penalization addresses it, we will train three neural networks to solve Burger's equations for a Riemann problem with hyperparameters
$\epsilon_{ent}\in\{0,0.01,1.0\}$, and $\epsilon_AV=\epsilon_{TVD}=0$. 
For Burger's the entropy pair used is 
\begin{equation}
    \eta = u^2; \quad q = \frac{2}{3} u^3.
\end{equation}
The initial network weights and biases are selected so that the initial guess is equal to a rarefaction shock solution
\begin{equation} \label{eq:rareshock}
    u = \left\{ \begin{matrix} 
            u_L = 1 & \mathrm{if} \ x <  \frac{1}{2} t, \\
            u_R = 0 & \mathrm{if} \ x \geq \frac{1}{2} t.
    \end{matrix} \right.
\end{equation}
We perform this initialization by choosing a set of collocation points and using gradient descent to minimize the $\ell_2$ norm between the neural network $u$ and Eqn.~\ref{eq:rareshock} at the collocation points. Finally, each
network is trained using cvPINNs with a different scaling of the entropy penalization. Table~\ref{table:ent_param} lists the hyperparameters for the minimization procedures for initialization and cvPINNs. 

In Figure~\ref{fig:ent} we examine the effects of the entropy penalization on cvPINNs.
By construction the rarefaction shock solution satisfies Eqn.~\ref{eq:pde} for all time. When $\epsilon_{ent}=0$ is used (as in the left plot), the optimal neural network has already been obtained by the initial guess. As a result 
we find that training yields the nonphysical rarefaction solution, yet it violates Eqn.~\ref{eq:pdie} (shown in the second row of the image). However, for $\epsilon_{ent}=0.01$ or $1.0$, the entropy solution is recovered and no longer violates Eqn.~\ref{eq:eres}. Moreover, the solution is little changed between the two non-zero penalties. Thus, in the remainder we use $\epsilon_{ent} = 1$.

\subsection{TVD condition} \label{sec:tvd}

\begin{figure}[htpb]
    \centering
    \includegraphics[width=\linewidth]{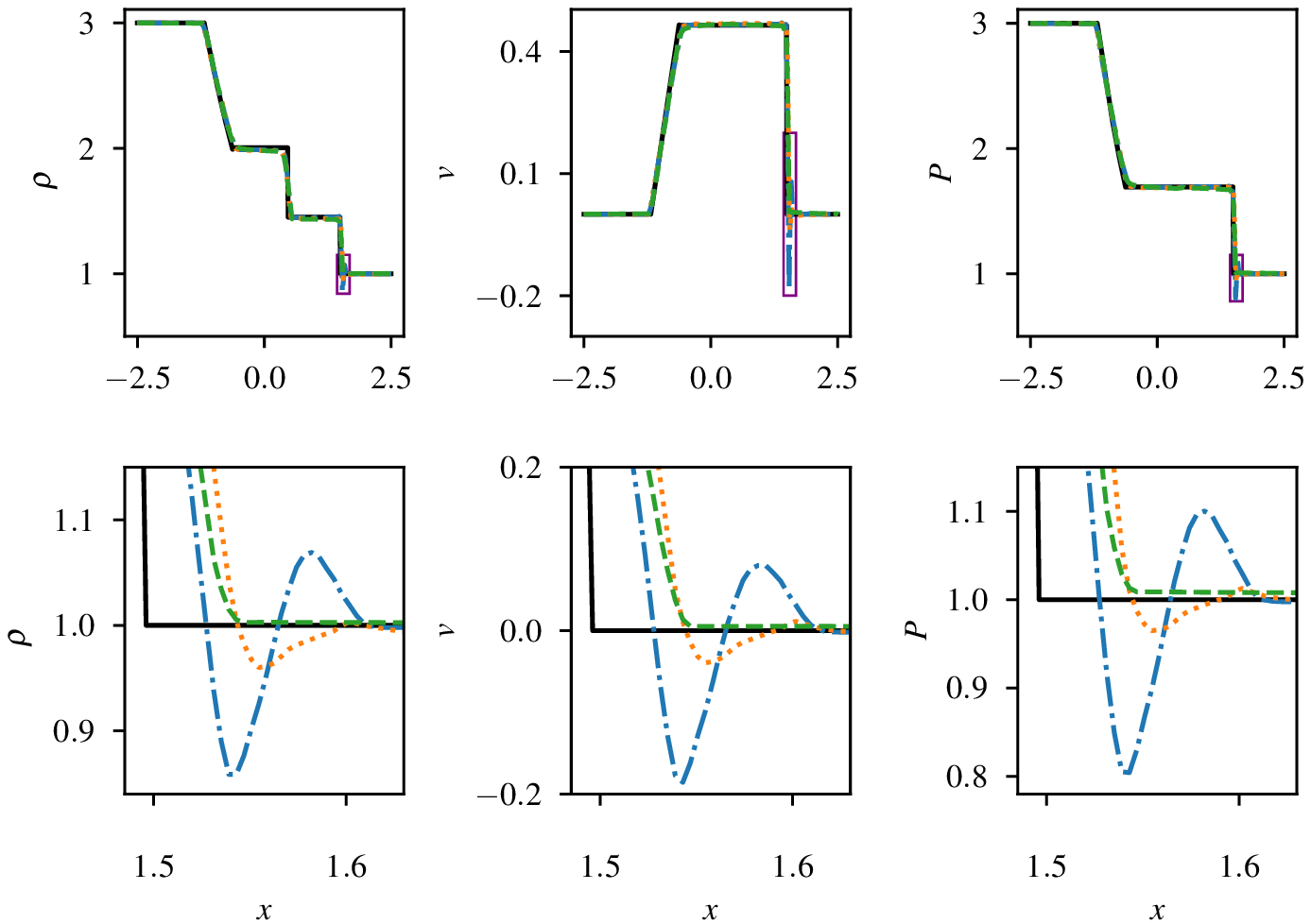}
    \caption{Effect of TVD regularization on cvPINNs solution for Sod shock problem. (\textit{top row}) Density (\textit{first column}), velocity (\textit{middle column}), and pressure (\textit{last column}) profiles at time $t=1$ for TVD regularization weights, $\epsilon_{TVD}=0$ \sym{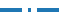}, $\epsilon_{TVD}=0.01$ \sym{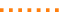}, and $\epsilon_{TVD}=1$ \sym{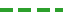}, and analytical solution \sym{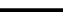}. (\textit{bottom row}) Zoomed in profiles at shock front. The profiles shown are the solutions with minimal loss over 10 trials. Here the
    reduction in oscillations is apparent for increasing values of $\epsilon_{TVD}$.}
    \label{fig:tvd}
\end{figure}


To explore the effect of the TVD penalization, cvPINNs is used to solve a Riemann problem using Euler's equations (Eqn.~\ref{eq:euler}) assuming a perfect gas EOS (Eqn.~\ref{eq:eosPgas} with $\gamma=1.4$). 
Figure~\ref{fig:tvd} shows the cvPINNs solution to the Sod problem, 
\begin{equation}\label{eq:riemann}
    \begin{bmatrix} \rho \\ v \\ p \end{bmatrix}_{t=0} = \begin{bmatrix} 3 \\ 0 \\ 3 \end{bmatrix}
    \quad \mathrm{if} \ x<0; \qquad 
    \begin{bmatrix} \rho \\ v \\ p \end{bmatrix}_{t=0} = \begin{bmatrix} 1 \\ 0 \\ 1 \end{bmatrix}
    \quad \mathrm{if} \ x\geq 0.
\end{equation}
using $\epsilon_{TVD}=0$, $0.01$ and $1.0$. 
Without the TVD penalization, using $\epsilon_{TVD}=0$, the plots demonstrate that cvPINNs solutions suffer from substantial oscillations, particularly at the shock front. However, the solutions arising from the non-zero values of $\epsilon_{TVD}$ show physical oscillations are incrementally reduced by increasing $\epsilon_{TVD}$. Here, the solution is essentially non-oscillatory for $\epsilon_{TVD}=1$ and does not appear to suffer from further deleterious dissipation effects. Thus, for the remainder of the paper, we use $\epsilon_{TVD} = 1$ when applying the TVD penalization. The hyperparameters used for cvPINNs are shown in Table~\ref{table:tvd_param}.

\subsection{Artificial Viscosity} \label{sec:av}

\begin{figure}[htpb]
    \centering
    \includegraphics[width=\linewidth]{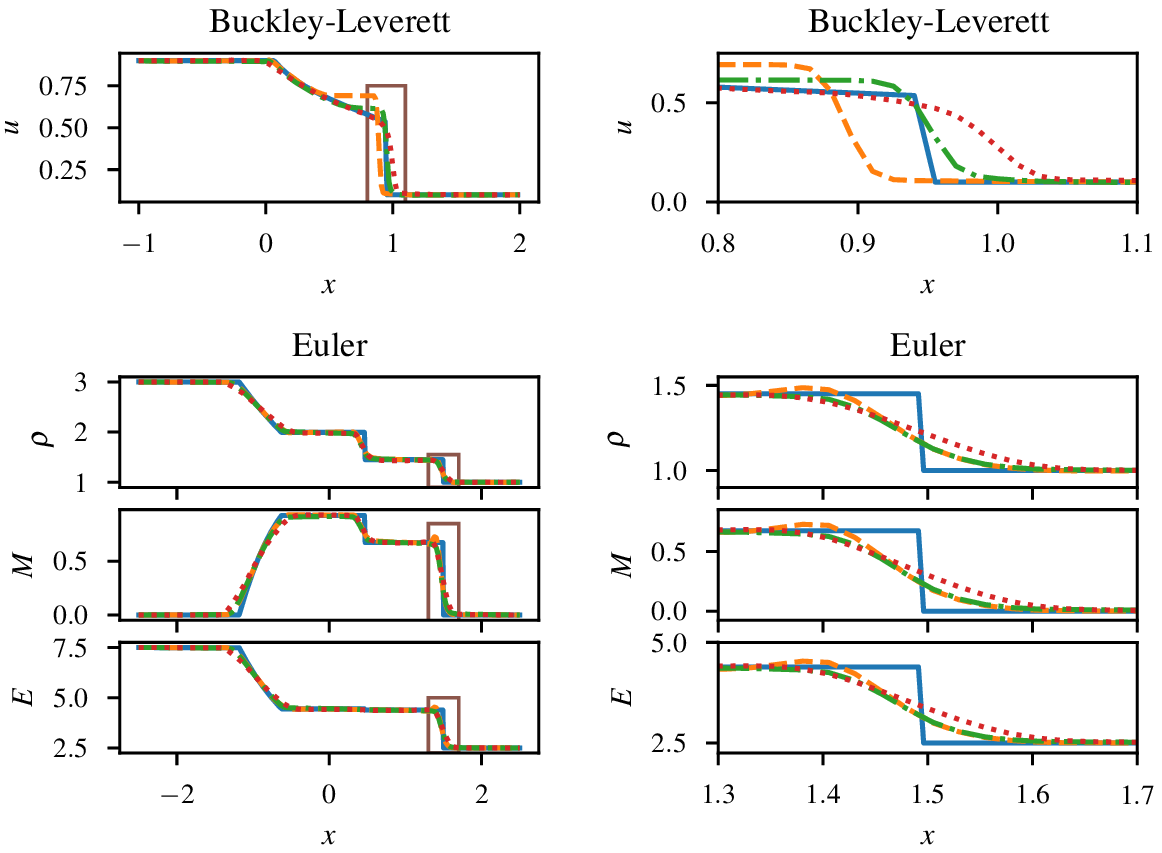}
    \caption{Entropy+TVD penalization and viscous penalization results in similar quality solutions for Buckley-Leverett and Euler Riemann problems. True solution \sym{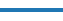}, no regularization  \sym{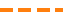}, entropy+TVD regularization  \sym{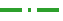}, viscous regularization \sym{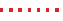}.}
    \label{fig:forward}
\end{figure}

Section~\ref{sec:cvpinns} introduced an artificial viscosity penalization to enforce both the entropy condition and minimize oscillations (using $\epsilon_{AV}>0$ and $\epsilon_{TVD}=\epsilon_{ent}=0$ in Eqn~\ref{eq:pderesTVD}). This section compares artificial viscosity to cvPINNs using TVD and entropy penalization (e.g. $\epsilon_{TVD}=\epsilon_{ent}=1$ and $\epsilon_{AV}=0$ in Eqn~\ref{eq:pderesTVD}).
Riemann problems for the Euler equations with the perfect gas EOS (see Eqn.~\ref{eq:euler}) and the the Buckley-Leverett equation,
%
\begin{equation}
    \partial_t  u + \partial_x \left( \frac{u^2}{u^2 + \frac{1}{2}(1-u)^2} \right) = 0,
\end{equation}
with the entropy pair \cite{Kivva2020},
\begin{equation}
    \begin{split}
        &\eta = \frac{1}{2}u^2, \\
        &q =
\frac{2}{9} \left(\frac{u-2}{1 - 2u + 3u^2} 
    + \frac{1}{\sqrt{2}}\arctan\left(\frac{-1 + 3u}{\sqrt{2}}\right)
                      - \log (1 - 2u + 3u^2)
                    \right).
    \end{split}
\end{equation}
are solved and compared to analytic solutions. The methods for finding the analytical entropy solutions for these equations are available in \cite{LeVeque2002}. As the Buckley-Leverett equation has a non-convex flux the solution to the Riemann problems consists of a composite wave, consisting of a rarefaction connected to a shock.

In Figure~\ref{fig:forward}, we show solutions to unpenalized cvPINNs, cvPINNs using TVD and entropy penalization, artificial viscosity penalization, and the analytical solutions to Riemann problems. 
For Buckley-Leverett without penalization, the cvPINNs solution finds an incorrect solution with the wrong shock speed and an overshoot behind the shock. Adding the viscous penalty leads cvPINNs to a better, albeit more dissipative, solution. Using the entropy and TVD penalty results in a solution that recovers the correct shock speed, but retains a slight overshoot.

For the Euler equations, omitting penalization does find the correct entropy solution, however it has substantive overshoots at the shock. Adding  artificial viscosity dampens the oscillations, resulting in a more dissipated solution. Using the entropy and TVD penalty instead eliminates the oscillation and is less dissipative than the viscous penalization.
The hyperparameters used for minimization are available in Table~\ref{tab:gen}. For the remainder of the paper entropy and TVD penalization are used for cvPINNs solutions.

\subsection{Unstructured mesh} \label{sec:unstructured}

\begin{figure}[htpb]
    \centering
    \includegraphics[width=\linewidth]{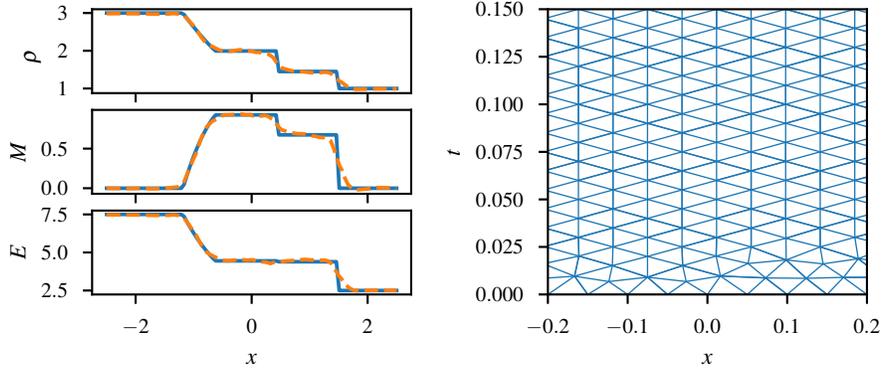}
    \caption{cvPINNs is compatible with a more general, unstructured mesh. (\textit{right}) Analytical solution to Eqn.~\ref{eq:riemann} \sym{C0,l,n,1} and solution using cvPINNs with entropy and TVD regularizations \sym{C1,d,n,1}. (\textit{right}) Section of the triangular mesh used with cvPINNs.}
    \label{fig:unstructured}
\end{figure}

Excluding this section, all of the tests consider Cartesian space-time meshes. However, cvPINNs is applicable on general polygonal meshes. To demonstrate this, we consider a Sod shock tube is solved on an unstructured, triangular space-time mesh, generated using pyGMSH~\cite{pygmsh} with characteristic lengths, $\Delta x/L = \Delta t/T = 100$. We apply identical parameters to the Cartesian case (see  Table~\ref{tab:gen}), and compare the resulting solution on a regular Cartesian grid of comparable size, $100\times 100$, demonstrating similar solution on structured/unstructured meshes (Figure~\ref{fig:unstructured}).

\subsection{Comparison of cvPINNs vs. PINNs}

\begin{figure}[htpb]
    \centering
    \includegraphics[width=\linewidth]{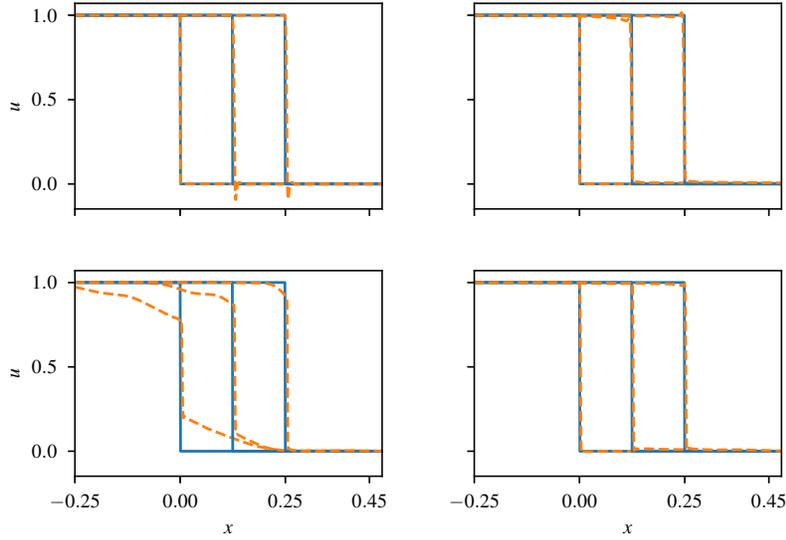}
    \caption{Comparison of cvPINNs (\textit{top row}) and traditional PINNs (\textit{bottom row}) for the Burger's shock problem demonstrate the sensitivity of PINNs to choice of IC/BC weighting. Analytical solution \sym{C0,l,n,1} and PINNs/cvPINNs solutions  \sym{C1,d,n,1} are shown for times, $t=0, 0.25, 0.5$.  (\textit{top left}) cvPINNs with no entropy or TVD regularization. (\textit{top right}) cvPINNs with entropy and TVD regularization weights set to 1. (\textit{bottom left}) PINNs with weights for IC and BC set to 1. (\textit{bottom right}) PINNs with weights for IC and BC set to 90. The profiles shown are the solutions with minimal loss over 10 trials.}
    \label{fig:PINNs}
\end{figure}


To highlight differences between cvPINNs and PINNs we consider the Riemann problem for the Burgers equation on the domain $[0,T]\times[-L/2,L/2]$,
\begin{equation}\label{eq:burgers}
    \begin{split}
    &\partial_t u + \frac{1}{2}\partial_x u^2 = 0, \\
    &u(0,x) = u_L \quad \mathrm{for} \quad x<0, \\
    &u(0,x) = u_R \quad \mathrm{for} \quad x\geq 0.
    \end{split}
\end{equation}
For $u_L>u_R$, the weak solution consists of a right moving shock with speed $\frac{u_L+u_R}{2}$ separating the left and right states. As is typically for nonlinear hyperbolic PDEs, the solution contains a discontinuity.

We compare PINNs and cvPINNs for solving the Burgers equation in Figure \ref{fig:PINNs}. Varying the penalty parameters weighting BC/IC demonstrates a large sensitivity in the solution. The cvPINNs solution correctly identifies the shock front, however incorporation of the TVD penalty is necessary to avoid oscillations at the shock front. 

Previous works have demonstrated some success with PINNs for nonlinear hyperbolic PDEs, although they often introduce regularization with sensitive parameters. In addition to the IC and BC penalty weight, \cite{Mao2020,Jagtap2020} successfully applied variations of PINNs to hyperbolic PDEs but require clustering collocation points around shocks. This necessitates \textit{a priori} knowledge of the shock location. Reference~\cite{Fuks2020} applied PINNs to the Buckley-Leverett equation and found that artificial viscosity was necessary to recover the entropy solution. However, the quality of the solution is sensitive to the choice of viscosity parameter. In constrast, cvPINNs with the TVD and entropy regularization appears to be robust with respect to the choice of parameter.

\section{Inverse problems: equations of state for shock physics}

In this section, we infer EOS's from solutions to the Euler equations using cvPINNs as discussed in Section~\ref{sec:inv}. In Section~\ref{sec:dsmc} we verify the cvPINNs method for inverse problems by extracting EOS's from well studied systems and comparing them to known EOS's. In Section~\ref{sec:lammps} we extract EOS's for shock hydrodynamics of copper. 

\subsection{Sod's shock tube from direct simulation Monte Carlo data} \label{sec:dsmc}

In this section, we apply cvPINNs to extract EOS's from simulations of compressible gas dynamics. We generate direct simulation Monte Carlo (DSMC) data sets of the Sod shock tube configuration \cite{Sod1978} for Argon under continuum flow conditions. DSMC solves the Boltzmann equation stochastically by using probabilistic descriptions of molecular behavior \cite{Bird1994}.  Since the solution of the Boltzmann equation is valid in continuum flows, DSMC can simulate continuum flows at low Knudsen number ($Kn < 1$) \cite{Gallis2015,Gallis2016,Gallis2018}.
For this study, we use SPARTA 
\cite{Plimpton2019} to generate DSMC data. The domain used to simulate each shock tube is $1\mathrm{mm} \times 50 \mathrm{mm}$ tiled by  $8\times 10^4$ cells. Initially, a hundred DSMC simulators are allocated to each cell. A variable-soft-sphere (VSS) model for collision dynamics between DSMC simulators is used. Finally, the time interval between collisions is fixed to $100 \mathrm{ps}$.

From the DSMC data, we extract 1D density, momentum, and total energy profiles. This includes noise due to the stochastic nature of DSMC. To enhance the signal-to-noise ratio of the DSMC data, a box filter is applied to downsample the profiles by a factor of 300.

\begin{figure}[htpb]
    \centering
    \includegraphics[width=4in]{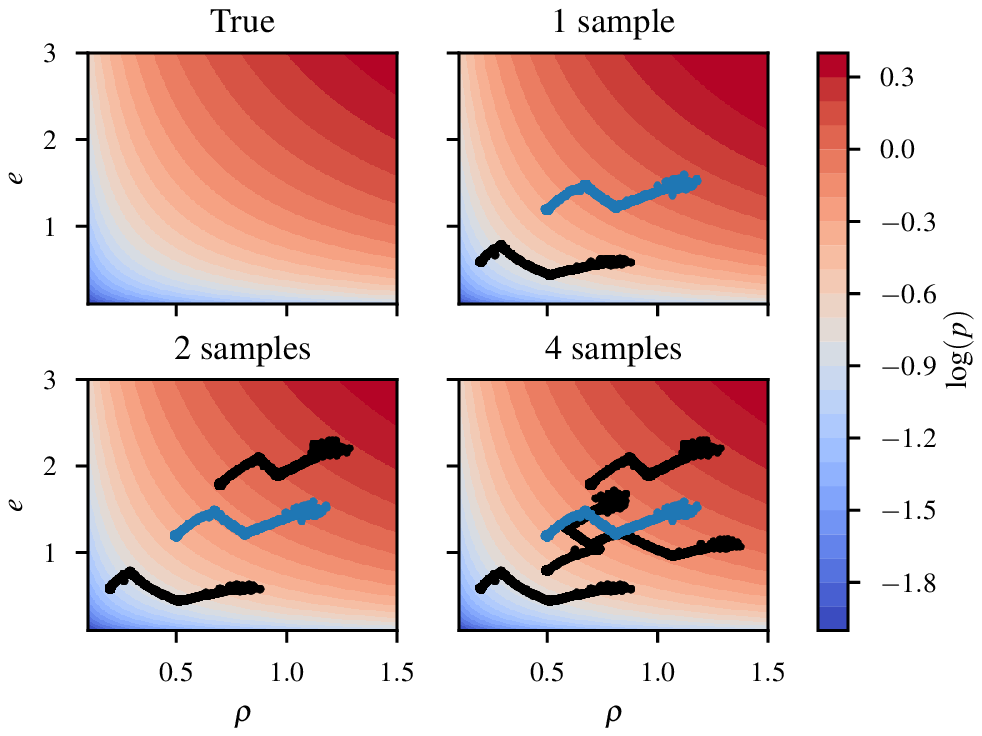}
    \caption{Pressure from regressed EOS from DSMC data using parameterized perfect gas. Perfect gas for argon (\textit{Top left}, $\gamma=5/3$). Learned EOS using data from one (\textit{Top right}, $\gamma=1.684$), two (\textit{Bottom left}, $\gamma=1.679$), and four (\textit{Bottom right}, $\gamma=1.677$) Riemann problems.  (\textit{Bottom left}). Data samples \sym{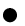} for regression, states for Riemann problem interpolated in Figure~\ref{fig:dsmc_interp} \sym{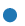}. This EOS parameterization makes strong assumptions but yields good generalizability.}
    \label{fig:dsmc_eos_ideal}
\end{figure}

\begin{figure}[htpb]
    \centering
    \includegraphics[width=4in]{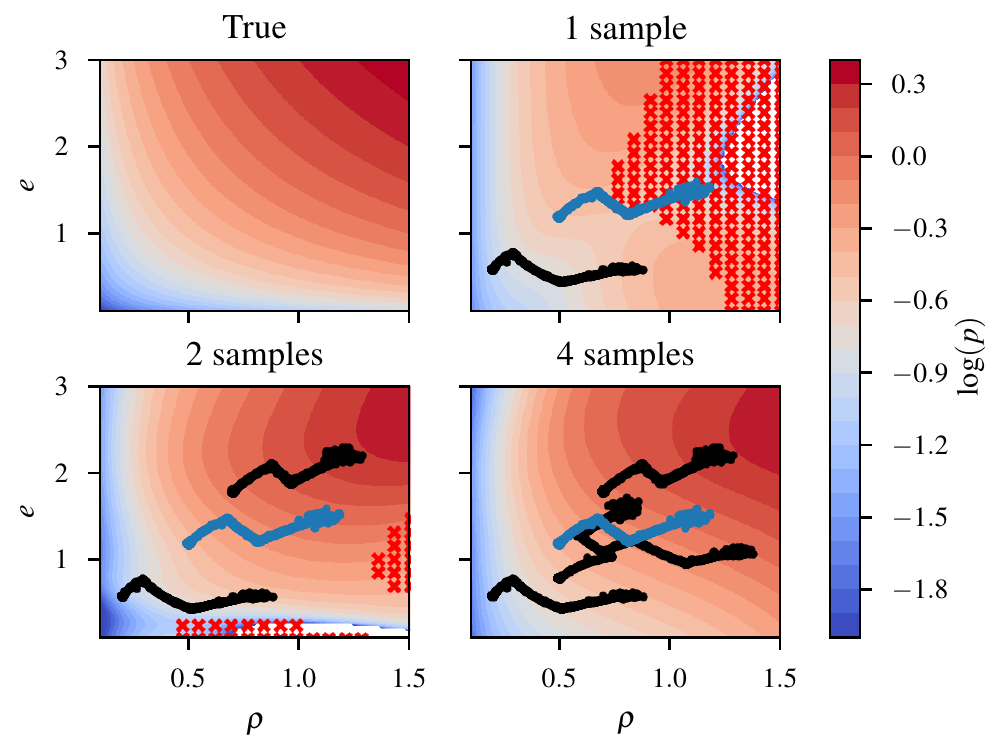}
    \caption{Pressure from regressed EOS from DSMC data using neural network. Perfect gas for argon (\textit{Top left}). Learned EOS using data from one (\textit{Top right}), two (\textit{Bottom left}), and four (\textit{Bottom right}) Riemann problems.  (\textit{Bottom left}). Data samples \sym{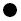} for regression, states for Riemann problem interpolated in Figure~\ref{fig:dsmc_interp} \sym{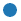}. Elliptic regions for $u=0$ \sym{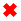}. This EOS parameterization makes few assumptions but yields lacks generalizability.}
    \label{fig:dsmc_eos_NN}
\end{figure}

\begin{figure}[htpb]
    \centering
    \includegraphics[width=4in]{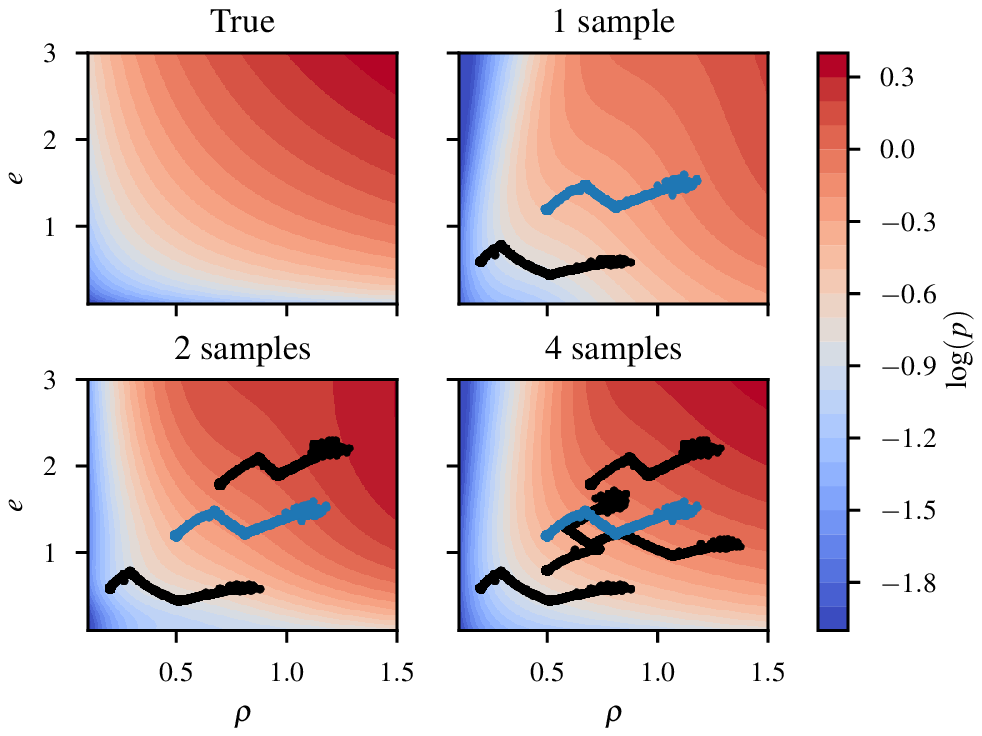}
    \caption{Pressure from regressed EOS from DSMC data using neural network with regularization. Perfect gas for argon (\textit{Top left}). Learned EOS using data from one (\textit{Top right}), two (\textit{Bottom left}), and four (\textit{Bottom right}) Riemann problems.  (\textit{Bottom left}). Data samples  \sym{k,n,o,4} for regression, states for Riemann problem interpolated in Figure~\ref{fig:dsmc_interp} \sym{C0,n,o,4}. Elliptic regions for $u=0$ \sym{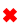}. This EOS parameterization balances generalizability and assumptions on the EOS form.}
    \label{fig:dsmc_eos_NN_reg}
\end{figure}

\begin{figure}[htpb]
    \centering
    \includegraphics[width=4.5in]{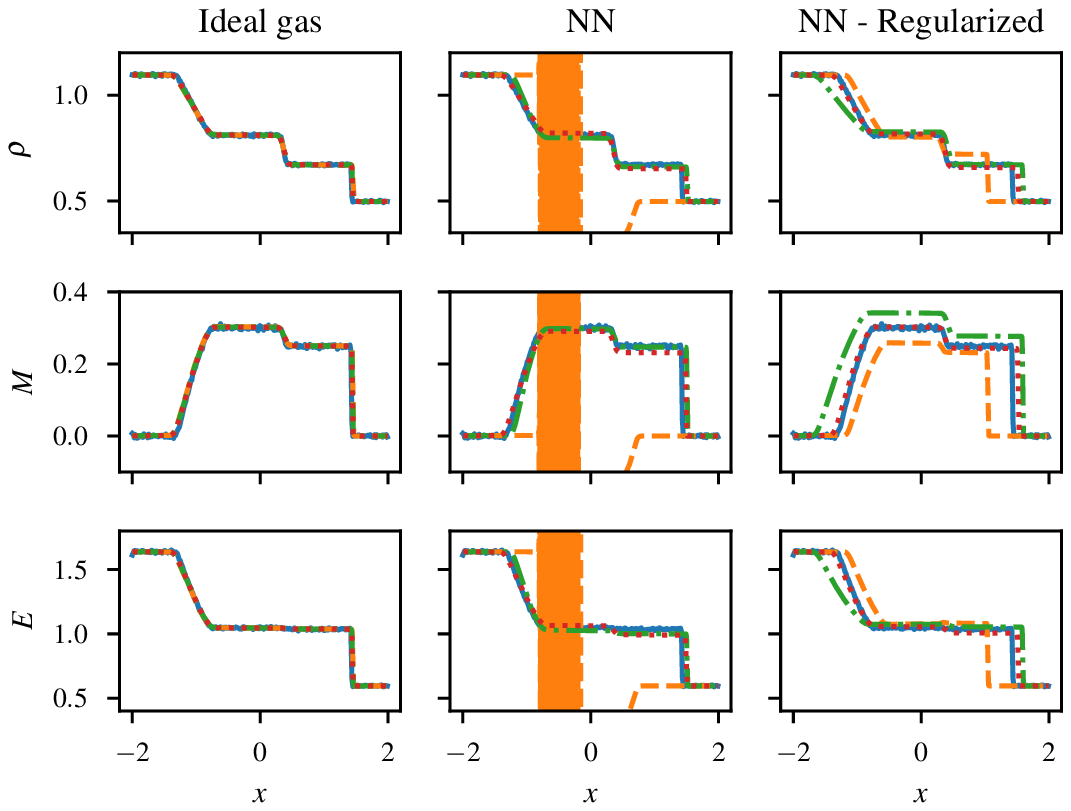}
    \caption{Solutions to test case Riemann problem using fitted EOS. Parameterized perfect gas (\textit{left column}), neural network (\textit{middle column}), and regularized neural network (\textit{right column}). Finite difference solution for fits using data from one \sym{C1,d,n,1}, two \sym{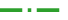}, and four \sym{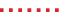} Riemann problems. DSMC solution \sym{C0,l,n,1}. Due to loss of hyperbolicity, the NN EOS parameterization can result in an unstable discretized PDE. The regularized parameterization always gives a stable discretized PDE. The test case accuracy improves for the NN parameterizations as more training data is used.}
    \label{fig:dsmc_interp}
\end{figure}

Five Riemann problems are simulated using SPARTA with varying pressure and density jumps. The data is normalized to the length of the domain, $L=1 \textrm{mm}$, the total time of the simulation, $1\mu$s, and the density of argon at standard temperature and pressure, $1.449$kg/m$^3$. For the regime considered, the Euler Riemann problem with the perfect gas EOS ($\gamma = 5/3$) provides a good model for the dynamics. We recover EOS's from the data and compare them to the perfect gas EOS with $\gamma=5/3$. The left and right states are available in Table~\ref{tab:dsmc}. We use the first four cases of this data set as a training set and the last as a test set. The hyperparameters used to train the EOS's is available in Table~\ref{tab:inv}. 

To assess the quality of EOS obtained, we take the trained EOS and implement it in a traditional finite-difference (FD) based Euler solver (for details, see \ref{app:fd}). This allows an assessment of the accuracy of the equation of state when deployed into a more traditional solver.

In Figure~\ref{fig:dsmc_eos_ideal}, we examine the recovered EOS's using the perfect gas EOS parameterization in Eqn.~\ref{eq:eosPgas}. We are able to recover the EOS with low error with only a single sample of DSMC shock tube data. In the first column of Figure~\ref{fig:dsmc_interp}, the resulting FD simulation applied to unseen initial conditions and finds good agreement with DSMC data, regardless of the number of shock tube DSMC simulations used to learn the EOS. 

Figure~\ref{fig:dsmc_eos_NN} shows the recovered EOS's using the neural network parameterization without thermodynamic regularization, Eqn.~\ref{eq:resIR}. Here, with only a single DSMC shock tube simulation, we find an EOS that fits the data well, but generalizes poorly. We also find imaginary eigenvalues for the flux Jacobian for the states, $(\rho,e)$, in the figure with $u=0$. This suggests that the PDE is elliptic for those states and the initial boundary value problem (IBVP) is ill-posed. In the second column of Figure~\ref{fig:dsmc_interp}, we attempt to perform a FD simulation using the learned EOS on the test shock tube case with a state trajectory passing through the elliptic region. We find the FD solution to be unstable. However, as we add more DSMC shock tube simulations to out training, we are able to find EOS's that give hyperbolic PDEs for larger regions of state space and we observe the FD simulations of the test case to match well with the DSMC data as shown in Figure~\ref{fig:dsmc_interp}.

Figure~\ref{fig:dsmc_eos_NN_reg} shows the recovered EOS's using the neural network parameterization with thermodynamic regularization, Eqn.~\ref{eq:resIR}. With only one DMSC shock tube case in the training set, we find that the learned EOS poorly matches the reference perfect gas EOS. However, adding additional training data improves the accuracy of the learned EOS's. In Figure~\ref{fig:dsmc_interp}, we observe that the FD simulation with the learned EOS, although stable, matches poorly with the DMSC simulation when only one DSMC shock tube case was used for training. As more training data is used, the FD solution with the learned EOS matches better the DSMC data.

In summary, for the cases considered, we observe that if the EOS model form is known \textit{a priori}, only one solution is necessary to obtain a good fit. On the other end of the spectrum, if a completely black-box model form is used, the inverse problem provides an EOS that cannot be incorporated into a traditional Euler solver unless large amounts of training data are used. By incorporating thermodynamic consistency constraints however, we obtain a stable EOS even in small data limits.

\subsection{Shock hydrodynamics of copper from molecular dynamics data}\label{sec:lammps}

\begin{figure}[htpb]
    \centering
    \includegraphics[width=\linewidth]{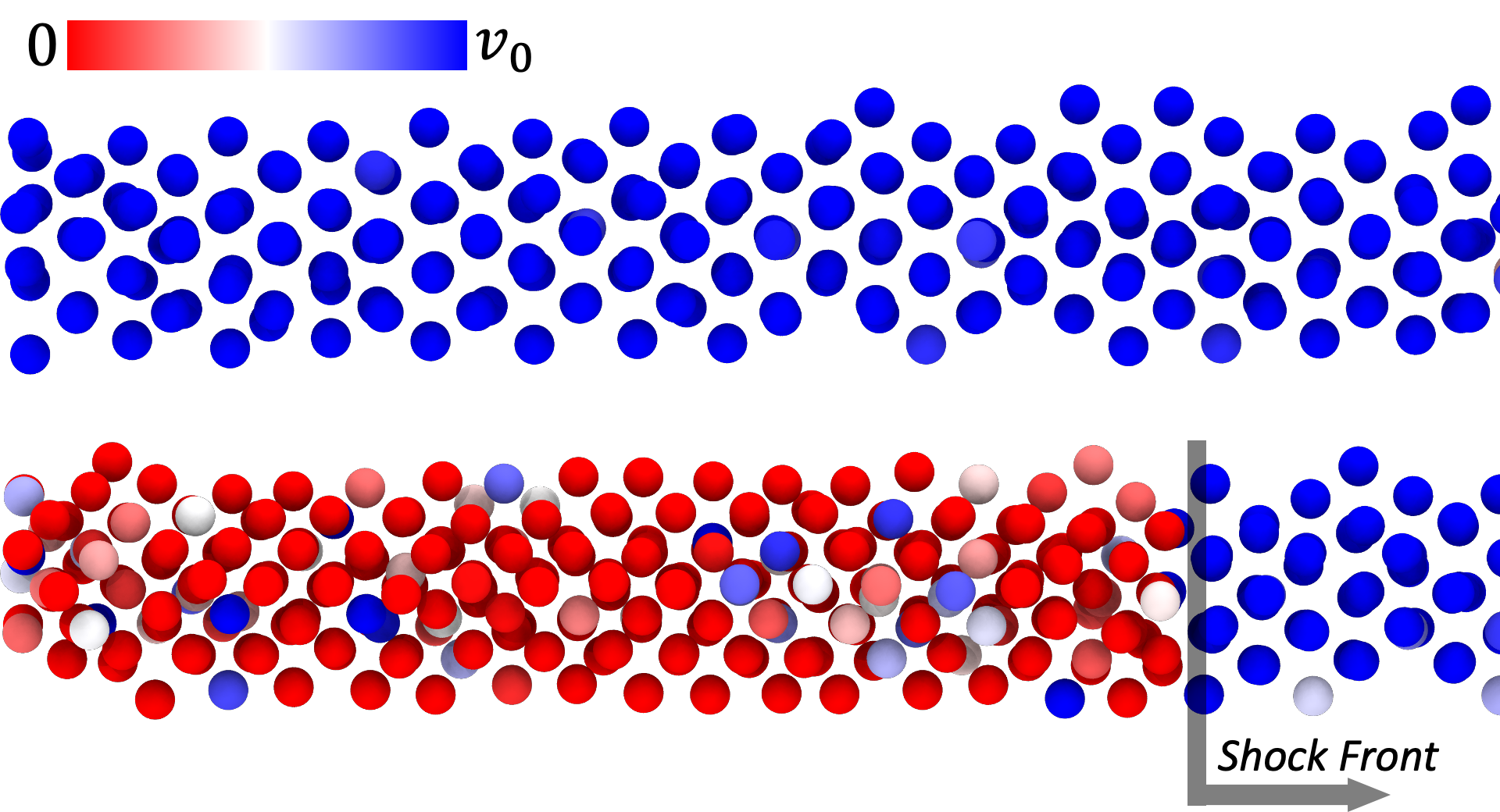}
    \caption{Visualization of copper impact LAMMPS molecular dynamics simulations rendered in Ovito \cite{ovito}. (\textit{top}) Copper bar at $t=0$. (\textit{bottom}) Copper bar shortly after impact. A right moving shock is observed via the velocity jump and the compression of copper atoms, $v_0=-2$km/s and $T_0=1000$K. }
    \label{fig:cu_vis}
\end{figure}

\begin{figure}[htpb]
    \centering
    \includegraphics[width=\linewidth]{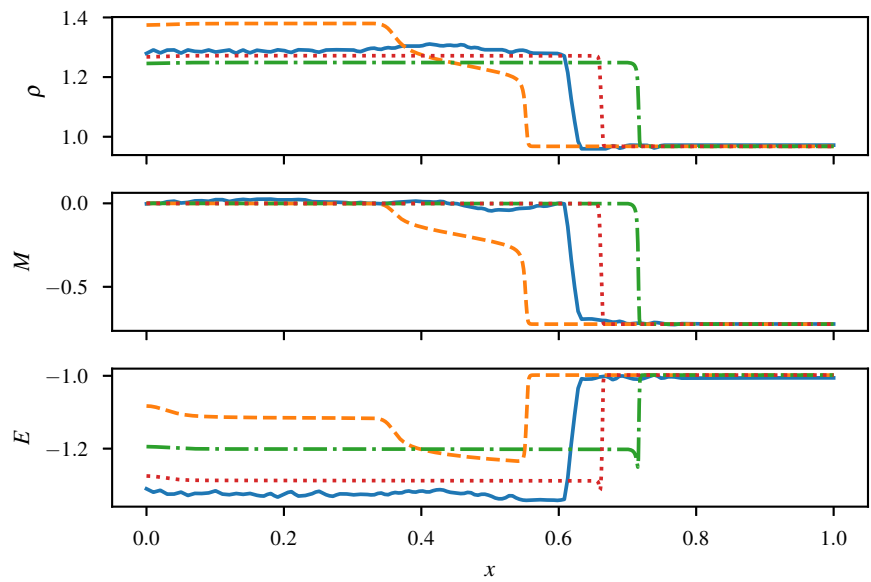}
    \caption{FD difference solutions for test copper impact case using EOS learned from training set with two cases \sym{C1,d,n,1}, four cases \sym{C2,od,n,1}, and eight cases \sym{C3,o,n,1}. DSMC solution for test case (\textit{blue}) is shown for comparison. The accuracy of the test case improves as more training data is used.}
    \label{fig:cu}
\end{figure}

\begin{figure}[htpb]
    \centering
    \includegraphics[width=3in]{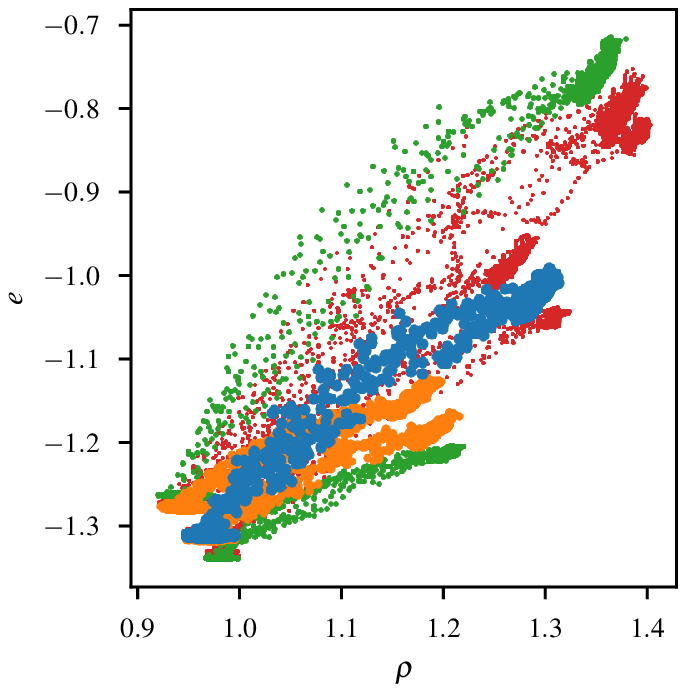}
    \caption{Projection of training data states onto ($\rho,e$) plane for two case training set \sym{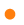}, four case training set (\syms{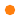},\syms{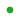}), eight case training set (\syms{C1,n,o,3},\syms{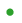},\syms{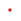}) and projection of test data \sym{C0,n,o,4}. As more training data is included, the training data's $(\rho,e)$ states surround the test data's and the test case switches from extrapolation to interpolation.}
    \label{fig:cu_samp}
\end{figure}

\begin{figure}[htpb]
    \centering
    \includegraphics[width=\linewidth]{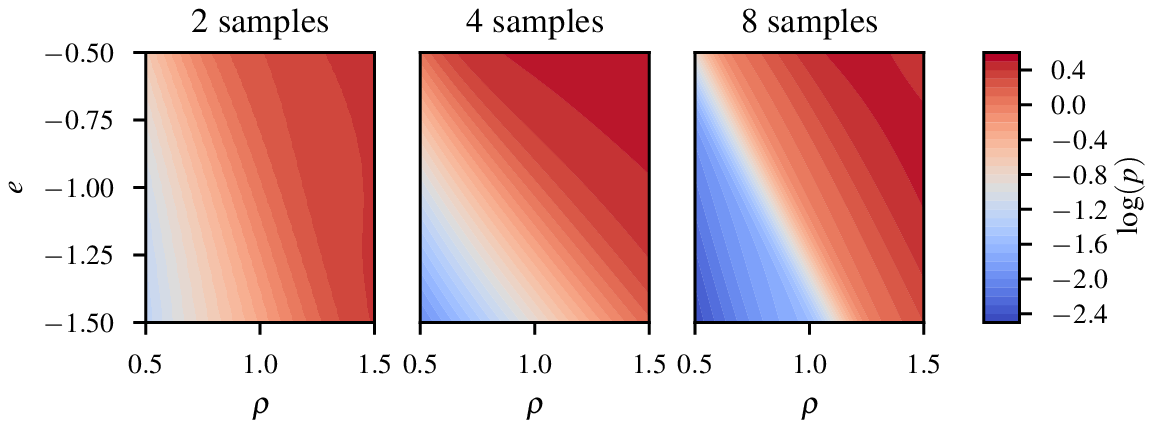}
    \caption{Pressure from the recovered EOS's using the 2 sample, 4 sample, and 8 sample training sets. Within the convex hull of the training points (See \ref{fig:cu_samp}) there is good agreement with the training data.}
    \label{fig:cu_eos}
\end{figure}

In this section, we apply the thermodynamic regularized cvPINNs method for inverse problems to infer the EOS for shock hydrodynamics of copper. We generate training and test data from LAMMPS molecular dynamics simulations of single crystal copper in a  reverse-ballistic impact experiment at various temperatures and shock impact velocities with the material described by an embedded atom interatomic potential \cite{plimpton1995fast, wood2018multiscale, mishin2001structural}. 
Hugoniot states of this interatomic potential have been previously shown to agree well with experimental data sets\cite{bringa2004atomistic}.

We performed nine LAMMPS simulations of the copper impact problem at varying temperatures and impact velocities and normalize the data using the length of the simulation domain, $1.4388\mu$m, and the speed of sound and density of copper at standard temperature and pressure, $3933$m/s and $8960$kg/m$^3$. The initial temperatures and impact velocities are given in Table~\ref{tab:lammps}. 
All shocks were aligned with the $[001]$ crystallographic direction which results in a single over driven elastic-plastic wave structure \cite{bringa2004atomistic}.
Each of the simulation geometries were constructed by replicating the FCC unit cell (2000:1 aspect ratio along shock direction) and equilibrating at the desired temperature while maintaining a constant pressure of one atmosphere with a Nose-Hoover barostat \cite{tuckerman2006liouville, parrinello1981polymorphic}.
As the shock direction is non-periodic, a stable starting pressure and temperature were achieved through serial relaxation steps of decreasing barostat dampening constants, an example input script is supplied in the supplemental material. 
Higher shock pressures and/or initial temperatures also result in melting.
An example visualization of given in Figure \ref{fig:cu_vis} wherein atoms are colored by the velocity in the shock direction, the simulated material will impact a stationary wall on the left-most boundary leaving shocked material stationary at a higher $\rho$ and $e$. 
We separate these simulations into a training set, with eight cases, and a test set, with just one case. We further separate the training set into three training sets of varying size, (i) with two shock conditions, (ii) four shock conditions, and (iii) with all eight. EOS's are learned with the thermodynamic regularized neural network parameterization with these training sets. The hyperparameters used for training are available in Table~\ref{tab:inv}. 

From these MD simulations, we extract 1D profile of density, momentum, and total energy. The copper impact simulations can be modeled with the 1D Euler equations with the IBVP,
\begin{equation}
    \begin{aligned}
        \begin{bmatrix} \rho \\ M \\ E \end{bmatrix}_{t=0} = 
        \begin{bmatrix} \rho_0 \\ M_0 \\ E_0 \end{bmatrix}; \quad  \begin{bmatrix} \partial_x \rho \\ M \\ \partial_x e \end{bmatrix}_{x=0} = 0; \quad
        \begin{bmatrix} \rho \\ M \\ E \end{bmatrix}_{x=L} = \begin{bmatrix} \rho_0 \\ M_0 \\ E_0 \end{bmatrix}
\end{aligned}
\end{equation}
with $M_0<0$. For cvPINNs, we use a finite difference approximation for the fluxes at the boundary, $x=0$, for the density and energy,
\begin{equation}
    \rho_{\frac12,n} = \rho_{1,n}; \quad e_{\frac12,n} = e_{1,n}.
\end{equation}

Figure~\ref{fig:cu} shows the FD solutions using the learned EOS for the test case problem and Figure~\ref{fig:cu_samp} shows a visualization of the training and test sets. The states in the first training set are far from the test states (see Figure \ref{fig:cu_samp}), as such we find that the FD solution for this EOS does not match well with the DSMC test solution. In fact, we observe a split shock wave instead of the single jump shock of the DSMC data. However, as we use the larger training sets that encompass a larger $\rho-e$ space to learn the EOS, we find improvements in the match between the FD solutions for the test problem with the learned EOS's and the DSMC test solution. Figure~\ref{fig:cu_eos} shows the pressures from the recovered EOS's from the three training sets.
The computational expense to generate the MD training data is modest, roughly 600 cpu$\cdot$hours, but has resulted in a general use EOS for a wide range of shock conditions.

\section{Conclusion}

We have presented a new physics-inform machine learning framework that incorporates a space-time control volume scheme to obtain significant improvements in accuracy and solution quality when solving conservation laws with neural networks. We have generated a number of new regularizers to introduce thermodynamic inductive biases and use these biases to solve inverse problems for equations of state over realistic, noisy data-sets corresponding to practical engineering problems.

The future of using deep learning architectures to solve numerical PDEs will require handling of the optimization error. While we have presented here a number of techniques to obtain qualitatively correct and physically meaningful solutions, the barrier in achieving convergence of error with respect to neural network size remains a major challenge to obtaining DNN solutions competitive with traditional finite element/volume methods. We refer the interested reader to some of our ongoing work in this area \cite{cyr2019robust}.

A particularly exciting development in the field of scientific machine learning is the area of operator regression, see for example \cite{lu2019deeponet,you2020data,patelaphysics}. We anticipate that the framework introduced here, together with the potential to introduce thermodynamically consistent biases, will allow cvPINNs to provide a fruitful area for learning more complex closures, such as those required for turbulence modeling, multiscale modeling, and obtaining closures for modeling non-equilibrium kinetics.

\section*{Acknowledgement}
Sandia National Laboratories is a multimission laboratory managed and operated by National Technology and Engineering Solutions of Sandia, LLC, a wholly owned subsidiary of Honeywell International, Inc., for the U.S. Department of Energy’s National Nuclear Security Administration under contract DE-NA0003525. This paper describes objective technical results and analysis. Any subjective views or opinions that might be expressed in the paper do not necessarily represent the views of the U.S. Department of Energy or the United States Government.

The work of R. Patel and N. Trask has also been supported by the U.S. Department of Energy, Office of Advanced Scientific Computing Research under the Collaboratory on Mathematics and Physics-Informed Learning Machines for Multiscale and Multiphysics Problems (PhILMs) project. The work of E. Cyr and N. Trask was supported by the U.S. Department of Energy, Office of Advanced Scientific Computing Research under the Early Career Research Program. 
SAND Number: SAND2020-13725 O.
The work of Myoungkyu Lee was supported by the US Department of Energy, Office of Basic Energy Sciences, Division of Chemical Sciences, Geosciences, and Biosciences. 

For generating data from DSMC, this research used resources of the Oak Ridge Leadership Computing Facility, which is a DOE Office of Science User Facility supported under Contract DE-AC05-00OR22725.

\bibliographystyle{elsarticle-num}
\bibliography{ref}{}

\begin{thebibliography}{10}
\expandafter\ifx\csname url\endcsname\relax
  \def\url#1{\texttt{#1}}\fi
\expandafter\ifx\csname urlprefix\endcsname\relax\def\urlprefix{URL }\fi
\expandafter\ifx\csname href\endcsname\relax
  \def\href#1#2{#2} \def\path#1{#1}\fi

\bibitem{weinan2018deep}
E.~Weinan, B.~Yu, The deep ritz method: a deep learning-based numerical
  algorithm for solving variational problems, Communications in Mathematics and
  Statistics 6~(1) (2018) 1--12.

\bibitem{he2018relu}
J.~He, L.~Li, J.~Xu, C.~Zheng, Relu deep neural networks and linear finite
  elements, arXiv preprint arXiv:1807.03973 (2018).

\bibitem{daubechies2019nonlinear}
I.~Daubechies, R.~DeVore, S.~Foucart, B.~Hanin, G.~Petrova, Nonlinear
  approximation and (deep) {ReLU} networks, arXiv preprint arXiv:1905.02199
  (2019).

\bibitem{yarotsky2017error}
D.~Yarotsky, Error bounds for approximations with deep {ReLU} networks, Neural
  Networks 94 (2017) 103--114.

\bibitem{yarotsky2018optimal}
D.~Yarotsky, Optimal approximation of continuous functions by very deep relu
  networks, arXiv preprint arXiv:1802.03620 (2018).

\bibitem{opschoor2019deep}
J.~A. Opschoor, P.~Petersen, C.~Schwab, Deep {ReLU} networks and high-order
  finite element methods, SAM, ETH Z{\"u}rich (2019).

\bibitem{bach2017breaking}
F.~Bach, Breaking the curse of dimensionality with convex neural networks, The
  Journal of Machine Learning Research 18~(1) (2017) 629--681.

\bibitem{bengio2000taking}
S.~Bengio, Y.~Bengio, Taking on the curse of dimensionality in joint
  distributions using neural networks, IEEE Transactions on Neural Networks
  11~(3) (2000) 550--557.

\bibitem{han2018solving}
J.~Han, A.~Jentzen, E.~Weinan, Solving high-dimensional partial differential
  equations using deep learning, Proceedings of the National Academy of
  Sciences 115~(34) (2018) 8505--8510.

\bibitem{wang2020understanding}
S.~Wang, Y.~Teng, P.~Perdikaris, Understanding and mitigating gradient
  pathologies in physics-informed neural networks, arXiv preprint
  arXiv:2001.04536 (2020).

\bibitem{beck2019full}
C.~Beck, A.~Jentzen, B.~Kuckuck, Full error analysis for the training of deep
  neural networks, arXiv preprint arXiv:1910.00121 (2019).

\bibitem{fokina2019growing}
D.~Fokina, I.~Oseledets, Growing axons: greedy learning of neural networks with
  application to function approximation, arXiv preprint arXiv:1910.12686
  (2019).

\bibitem{adcock2020gap}
B.~Adcock, N.~Dexter, The gap between theory and practice in function
  approximation with deep neural networks, arXiv preprint arXiv:2001.07523
  (2020).

\bibitem{lagaris1998artificial}
I.~E. Lagaris, A.~Likas, D.~I. Fotiadis, Artificial neural networks for solving
  ordinary and partial differential equations, IEEE transactions on neural
  networks 9~(5) (1998) 987--1000.

\bibitem{raissi2019physics}
M.~Raissi, P.~Perdikaris, G.~E. Karniadakis, Physics-informed neural networks:
  A deep learning framework for solving forward and inverse problems involving
  nonlinear partial differential equations, Journal of Computational Physics
  378 (2019) 686--707.

\bibitem{raissi2018deep}
M.~Raissi, Deep hidden physics models: Deep learning of nonlinear partial
  differential equations, The Journal of Machine Learning Research 19~(1)
  (2018) 932--955.

\bibitem{sun2020surrogate}
L.~Sun, H.~Gao, S.~Pan, J.-X. Wang, Surrogate modeling for fluid flows based on
  physics-constrained deep learning without simulation data, Computer Methods
  in Applied Mechanics and Engineering 361 (2020) 112732.

\bibitem{zhang2019quantifying}
D.~Zhang, L.~Lu, L.~Guo, G.~E. Karniadakis, Quantifying total uncertainty in
  physics-informed neural networks for solving forward and inverse stochastic
  problems, Journal of Computational Physics 397 (2019) 108850.

\bibitem{meng2020composite}
X.~Meng, G.~E. Karniadakis, A composite neural network that learns from
  multi-fidelity data: Application to function approximation and inverse pde
  problems, Journal of Computational Physics 401 (2020) 109020.

\bibitem{mao2020physics}
Z.~Mao, A.~D. Jagtap, G.~E. Karniadakis, Physics-informed neural networks for
  high-speed flows, Computer Methods in Applied Mechanics and Engineering 360
  (2020) 112789.

\bibitem{zhang2020learning}
D.~Zhang, L.~Guo, G.~E. Karniadakis, Learning in modal space: Solving
  time-dependent stochastic pdes using physics-informed neural networks, SIAM
  Journal on Scientific Computing 42~(2) (2020) A639--A665.

\bibitem{tensorflow2015-whitepaper}
M.~Abadi, A.~Agarwal, P.~Barham, E.~Brevdo, Z.~Chen, C.~Citro, G.~S. Corrado,
  A.~Davis, J.~Dean, M.~Devin, S.~Ghemawat, I.~Goodfellow, A.~Harp, G.~Irving,
  M.~Isard, Y.~Jia, R.~Jozefowicz, L.~Kaiser, M.~Kudlur, J.~Levenberg,
  D.~Man\'{e}, R.~Monga, S.~Moore, D.~Murray, C.~Olah, M.~Schuster, J.~Shlens,
  B.~Steiner, I.~Sutskever, K.~Talwar, P.~Tucker, V.~Vanhoucke, V.~Vasudevan,
  F.~Vi\'{e}gas, O.~Vinyals, P.~Warden, M.~Wattenberg, M.~Wicke, Y.~Yu,
  X.~Zheng, {TensorFlow}: Large-scale machine learning on heterogeneous
  systems, software available from tensorflow.org (2015).

\bibitem{NEURIPS2019_9015}
A.~Paszke, S.~Gross, F.~Massa, A.~Lerer, J.~Bradbury, G.~Chanan, T.~Killeen,
  Z.~Lin, N.~Gimelshein, L.~Antiga, A.~Desmaison, A.~Kopf, E.~Yang, Z.~DeVito,
  M.~Raison, A.~Tejani, S.~Chilamkurthy, B.~Steiner, L.~Fang, J.~Bai,
  S.~Chintala, Pytorch: An imperative style, high-performance deep learning
  library, in: H.~Wallach, H.~Larochelle, A.~Beygelzimer, F.~d'~Alch\'{e}-Buc,
  E.~Fox, R.~Garnett (Eds.), Advances in Neural Information Processing Systems
  32, Curran Associates, Inc., 2019, pp. 8024--8035.

\bibitem{robinson2013fundamental}
A.~Robinson, R.~Berry, J.~Carpenter, B.~Debusschere, R.~R. Drake, A.~Mattsson,
  W.~J. Rider, Fundamental issues in the representation and propagation of
  uncertain equation of state information in shock hydrodynamics, Computers \&
  Fluids 83 (2013) 187--193.

\bibitem{carpenter2015automated}
J.~H. Carpenter, A.~C. Robinson, B.~Debusschere, A.~E. Wills, Automated
  generation of tabular equations of state with uncertainty information., Tech.
  rep., Sandia National Lab.(SNL-NM), Albuquerque, NM (United States);
  Sandia~… (2015).

\bibitem{jagtap2020conservative}
A.~D. Jagtap, E.~Kharazmi, G.~E. Karniadakis, Conservative physics-informed
  neural networks on discrete domains for conservation laws: Applications to
  forward and inverse problems, Computer Methods in Applied Mechanics and
  Engineering 365 (2020) 113028.

\bibitem{lax1973hyperbolic}
P.~D. Lax, Hyperbolic systems of conservation laws and the mathematical theory
  of shock waves, SIAM, 1973.

\bibitem{menikoff1989riemann}
R.~Menikoff, B.~J. Plohr, The riemann problem for fluid flow of real materials,
  Reviews of modern physics 61~(1) (1989) 75.

\bibitem{Meni1989}
R.~Menikoff, B.~J. Plohr,
  \href{https://doi-org.libproxy.unm.edu/10.1103/RevModPhys.61.75}{The
  {R}iemann problem for fluid flow of real materials}, Rev. Modern Phys. 61~(1)
  (1989) 75--130.
\newblock \href {https://doi.org/10.1103/RevModPhys.61.75}
  {\path{doi:10.1103/RevModPhys.61.75}}.
\newline\urlprefix\url{https://doi-org.libproxy.unm.edu/10.1103/RevModPhys.61.75}

\bibitem{Viscous2014}
J.-L. Guermond, B.~Popov,
  \href{https://doi-org.libproxy.unm.edu/10.1137/120903312}{Viscous
  regularization of the {E}uler equations and entropy principles}, SIAM J.
  Appl. Math. 74~(2) (2014) 284--305.
\newblock \href {https://doi.org/10.1137/120903312}
  {\path{doi:10.1137/120903312}}.
\newline\urlprefix\url{https://doi-org.libproxy.unm.edu/10.1137/120903312}

\bibitem{xiong2020roenets}
S.~Xiong, X.~He, Y.~Tong, R.~Liu, B.~Zhu, Roenets: Predicting discontinuity of
  hyperbolic systems from continuous data, arXiv preprint arXiv:2006.04180
  (2020).

\bibitem{tokareva2019machine}
S.~Tokareva, M.~J. Shashkov, A.~N. Skurikhin, Machine learning approach for the
  solution of the riemann problem in fluid dynamics, Tech. rep., Los Alamos
  National Lab.(LANL), Los Alamos, NM (United States) (2019).

\bibitem{magiera2020constraint}
J.~Magiera, D.~Ray, J.~S. Hesthaven, C.~Rohde, Constraint-aware neural networks
  for riemann problems, Journal of Computational Physics 409 (2020) 109345.

\bibitem{kharazmi2020hp}
E.~Kharazmi, Z.~Zhang, G.~E. Karniadakis, hp-vpinns: Variational
  physics-informed neural networks with domain decomposition, arXiv preprint
  arXiv:2003.05385 (2020).

\bibitem{jagtap2020extended}
A.~D. Jagtap, G.~E. Karniadakis, Extended physics-informed neural networks
  (xpinns): A generalized space-time domain decomposition based deep learning
  framework for nonlinear partial differential equations (2020).

\bibitem{Darfermos2000}
C.~M. Dafermos,
  \href{https://doi-org.libproxy.unm.edu/10.1007/3-540-29089-3_14}{Hyperbolic
  conservation laws in continuum physics}, Vol. 325 of Grundlehren der
  Mathematischen Wissenschaften [Fundamental Principles of Mathematical
  Sciences], Springer-Verlag, Berlin, 2000.
\newblock \href {https://doi.org/10.1007/3-540-29089-3_14}
  {\path{doi:10.1007/3-540-29089-3_14}}.
\newline\urlprefix\url{https://doi-org.libproxy.unm.edu/10.1007/3-540-29089-3_14}

\bibitem{Bianchi2005}
S.~Bianchini, A.~Bressan,
  \href{https://doi-org.libproxy.unm.edu/10.4007/annals.2005.161.223}{Vanishing
  viscosity solutions of nonlinear hyperbolic systems}, Ann. of Math. (2)
  161~(1) (2005) 223--342.
\newblock \href {https://doi.org/10.4007/annals.2005.161.223}
  {\path{doi:10.4007/annals.2005.161.223}}.
\newline\urlprefix\url{https://doi-org.libproxy.unm.edu/10.4007/annals.2005.161.223}

\bibitem{GodRav1996}
E.~Godlewski, P.-A. Raviart,
  \href{https://doi-org.libproxy.unm.edu/10.1007/978-1-4612-0713-9}{Numerical
  approximation of hyperbolic systems of conservation laws}, Vol. 118 of
  Applied Mathematical Sciences, Springer-Verlag, New York, 1996.
\newblock \href {https://doi.org/10.1007/978-1-4612-0713-9}
  {\path{doi:10.1007/978-1-4612-0713-9}}.
\newline\urlprefix\url{https://doi-org.libproxy.unm.edu/10.1007/978-1-4612-0713-9}

\bibitem{abadi2016tensorflow}
M.~Abadi, P.~Barham, J.~Chen, Z.~Chen, A.~Davis, J.~Dean, M.~Devin,
  S.~Ghemawat, G.~Irving, M.~Isard, et~al., Tensorflow: A system for
  large-scale machine learning, in: 12th $\{$USENIX$\}$ symposium on operating
  systems design and implementation ($\{$OSDI$\}$ 16), 2016, pp. 265--283.

\bibitem{moritz1978least}
H.~Moritz, Least-squares collocation, Reviews of geophysics 16~(3) (1978)
  421--430.

\bibitem{rummel1979least}
R.~Rummel, K.-P. Schwarz, M.~Gerstl, Least squares collocation and
  regularization, Bulletin Geodesique 53~(4) (1979) 343--361.

\bibitem{ling2005least}
L.~Ling, E.~J. Kansa, A least-squares preconditioner for radial basis functions
  collocation methods, Advances in Computational Mathematics 23~(1-2) (2005)
  31--54.

\bibitem{hu2007weighted}
H.~Hu, J.~Chen, W.~Hu, Weighted radial basis collocation method for boundary
  value problems, International journal for numerical methods in engineering
  69~(13) (2007) 2736--2757.

\bibitem{zhang2001least}
X.~Zhang, X.-H. Liu, K.-Z. Song, M.-W. Lu, Least-squares collocation meshless
  method, International Journal for Numerical Methods in Engineering 51~(9)
  (2001) 1089--1100.

\bibitem{cheng2010collocation}
H.~Cheng, A.~Sandu, Collocation least-squares polynomial chaos method, in:
  Proceedings of the 2010 Spring Simulation Multiconference, 2010, pp. 1--6.

\bibitem{BochevGLSbook}
P.~B. Bochev, M.~D. Gunzburger, Least-squares finite element methods, Vol. 166
  of Applied Mathematical Sciences, Springer, New York, 2009.
\newblock \href {https://doi.org/10.1007/b13382} {\path{doi:10.1007/b13382}}.

\bibitem{Jagtap2020}
A.~D. Jagtap, E.~Kharazmi, G.~E. Karniadakis,
  \href{https://www.sciencedirect.com/science/article/pii/S0045782520302127}{{Conservative
  physics-informed neural networks on discrete domains for conservation laws:
  Applications to forward and inverse problems}}, Computer Methods in Applied
  Mechanics and Engineering 365 (2020) 113028.
\newblock \href {https://doi.org/10.1016/J.CMA.2020.113028}
  {\path{doi:10.1016/J.CMA.2020.113028}}.
\newline\urlprefix\url{https://www.sciencedirect.com/science/article/pii/S0045782520302127}

\bibitem{Bochev2016}
P.~Bochev, M.~Gunzburger, Least-Squares Methods for Hyperbolic Problems, 2016.
\newblock \href {https://doi.org/10.1016/bs.hna.2016.07.002}
  {\path{doi:10.1016/bs.hna.2016.07.002}}.

\bibitem{L1Guermond}
J.-L. Guermond, F.~Marpeau, B.~Popov, A fast algorithm for solving first-order
  {PDE}s by {$L^1$}-minimization, Commun. Math. Sci. 6~(1) (2008) 199--216.

\bibitem{Guermond2008}
J.-L. Guermond, B.~Popov, {$L^1$}-minimization methods for {H}amilton-{J}acobi
  equations: the one-dimensional case, Numer. Math. 109~(2) (2008) 269--284.
\newblock \href {https://doi.org/10.1007/s00211-008-0142-1}
  {\path{doi:10.1007/s00211-008-0142-1}}.

\bibitem{Reisner2013}
J.~Reisner, J.~Serencsa, S.~Shkoller, {A space-time smooth artificial viscosity
  method for nonlinear conservation laws}, Journal of Computational Physics
  (2013).
\newblock \href {http://arxiv.org/abs/1204.0569} {\path{arXiv:1204.0569}},
  \href {https://doi.org/10.1016/j.jcp.2012.08.027}
  {\path{doi:10.1016/j.jcp.2012.08.027}}.

\bibitem{Rauch1986}
J.~Rauch, B{V} estimates fail for most quasilinear hyperbolic systems in
  dimensions greater than one, Comm. Math. Phys. 106~(3) (1986) 481--484.

\bibitem{Toro2000}
E.~F. Toro, S.~J. Billett,
  \href{https://doi-org.libproxy.unm.edu/10.1093/imanum/20.1.47}{Centred {TVD}
  schemes for hyperbolic conservation laws}, IMA J. Numer. Anal. 20~(1) (2000)
  47--79.
\newblock \href {https://doi.org/10.1093/imanum/20.1.47}
  {\path{doi:10.1093/imanum/20.1.47}}.
\newline\urlprefix\url{https://doi-org.libproxy.unm.edu/10.1093/imanum/20.1.47}

\bibitem{kaiser2018sparse}
E.~Kaiser, J.~N. Kutz, S.~L. Brunton, Sparse identification of nonlinear
  dynamics for model predictive control in the low-data limit, Proceedings of
  the Royal Society A 474~(2219) (2018) 20180335.

\bibitem{lu2019deeponet}
L.~Lu, P.~Jin, G.~E. Karniadakis, Deeponet: Learning nonlinear operators for
  identifying differential equations based on the universal approximation
  theorem of operators, arXiv preprint arXiv:1910.03193 (2019).

\bibitem{Harten1998}
A.~Harten, P.~D. Lax, C.~D. Levermore, W.~J. Morokoff,
  \href{https://doi-org.libproxy.unm.edu/10.1137/S0036142997316700}{Convex
  entropies and hyperbolicity for general {E}uler equations}, SIAM J. Numer.
  Anal. 35~(6) (1998) 2117--2127.
\newblock \href {https://doi.org/10.1137/S0036142997316700}
  {\path{doi:10.1137/S0036142997316700}}.
\newline\urlprefix\url{https://doi-org.libproxy.unm.edu/10.1137/S0036142997316700}

\bibitem{Kivva2020}
S.~Kivva, \href{http://arxiv.org/abs/2004.02258}{{Entropy stable flux
  correction for scalar hyperbolic conservation laws}} (apr 2020).
\newblock \href {http://arxiv.org/abs/2004.02258} {\path{arXiv:2004.02258}}.
\newline\urlprefix\url{http://arxiv.org/abs/2004.02258}

\bibitem{LeVeque2002}
R.~J. LeVeque, {Finite Volume Methods for Hyperbolic Problems}, 2002.
\newblock \href {https://doi.org/10.1017/cbo9780511791253}
  {\path{doi:10.1017/cbo9780511791253}}.

\bibitem{pygmsh}
N.~Schlömer, A.~Cervone, G.~D. McBain, tryfon mw, Nate, F.~Gokstorp, R.~van
  Staden, toothstone, J.~S. Dokken, D.~Kempf, J.~Sanchez, anzil, M.~Bussonnier,
  F.~Fu, ivanmultiwave, N.~Wagner, S.~Chen, tayebzaidi, T.~Maric, awa5114,
  Y.~Feng, \href{https://doi.org/10.5281/zenodo.4304309}{nschloe/pygmsh v7.1.5}
  (Dec. 2020).
\newblock \href {https://doi.org/10.5281/zenodo.4304309}
  {\path{doi:10.5281/zenodo.4304309}}.
\newline\urlprefix\url{https://doi.org/10.5281/zenodo.4304309}

\bibitem{Mao2020}
Z.~Mao, A.~D. Jagtap, G.~E. Karniadakis,
  \href{https://www.sciencedirect.com/science/article/pii/S0045782519306814}{{Physics-informed
  neural networks for high-speed flows}}, Computer Methods in Applied Mechanics
  and Engineering 360 (2020) 112789.
\newblock \href {https://doi.org/10.1016/J.CMA.2019.112789}
  {\path{doi:10.1016/J.CMA.2019.112789}}.
\newline\urlprefix\url{https://www.sciencedirect.com/science/article/pii/S0045782519306814}

\bibitem{Fuks2020}
O.~Fuks, H.~A. Tchelepi,
  \href{http://www.dl.begellhouse.com/journals/558048804a15188a,583c4e56625ba94e,415f83b5707fde65.html}{{LIMITATIONS
  OF PHYSICS INFORMED MACHINE LEARNING FOR NONLINEAR TWO-PHASE TRANSPORT IN
  POROUS MEDIA}}, Journal of Machine Learning for Modeling and Computing 1~(1)
  (2020).
\newblock \href {https://doi.org/10.1615/.2020033905}
  {\path{doi:10.1615/.2020033905}}.
\newline\urlprefix\url{http://www.dl.begellhouse.com/journals/558048804a15188a,583c4e56625ba94e,415f83b5707fde65.html}

\bibitem{Sod1978}
G.~A. Sod, A survey of several finite difference methods for systems of
  nonlinear hyperbolic conservation laws, Journal of Computational Physics 27
  (1978) 1--31.
\newblock \href {https://doi.org/10.1016/0021-9991(78)90023-2}
  {\path{doi:10.1016/0021-9991(78)90023-2}}.

\bibitem{Bird1994}
G.~A. Bird, Molecular Gas Dynamics and the Direct Simulation of Gas Flows,
  Oxford University Press, 1994.

\bibitem{Gallis2015}
M.~A. Gallis, T.~P. Koehler, J.~R. Torczynski, S.~J. Plimpton, Direct
  simulation monte carlo investigation of the richtmyer-meshkov instability,
  Physics of Fluids 27 (2015) 84105.
\newblock \href {https://doi.org/10.1063/1.4928338}
  {\path{doi:10.1063/1.4928338}}.

\bibitem{Gallis2016}
M.~A. Gallis, T.~P. Koehler, J.~R. Torczynski, S.~J. Plimpton, Direct
  simulation monte carlo investigation of the rayleigh-taylor instability,
  Physical Review Fluids 1 (2016) 43403.
\newblock \href {https://doi.org/10.1103/PhysRevFluids.1.043403}
  {\path{doi:10.1103/PhysRevFluids.1.043403}}.

\bibitem{Gallis2018}
M.~A. Gallis, J.~R. Torczynski, N.~P. Bitter, T.~P. Koehler, S.~J. Plimpton,
  G.~Papadakis, Gas-kinetic simulation of sustained turbulence in minimal
  couette flow, Physical Review Fluids 3 (2018) 71402.
\newblock \href {https://doi.org/10.1103/PhysRevFluids.3.071402}
  {\path{doi:10.1103/PhysRevFluids.3.071402}}.

\bibitem{Plimpton2019}
S.~J. Plimpton, S.~G. Moore, A.~Borner, A.~K. Stagg, T.~P. Koehler, J.~R.
  Torczynski, M.~A. Gallis, Direct simulation monte carlo on petaflop
  supercomputers and beyond, Physics of Fluids 31~(8) (2019) 086101.
\newblock \href {https://doi.org/10.1063/1.5108534}
  {\path{doi:10.1063/1.5108534}}.

\bibitem{ovito}
A.~Stukowski, {Visualization and analysis of atomistic simulation data with
  OVITO-the Open Visualization Tool}, {MODELLING AND SIMULATION IN MATERIALS
  SCIENCE AND ENGINEERING} {18}~({1}) ({JAN} {2010}).
\newblock \href {https://doi.org/{10.1088/0965-0393/18/1/015012}}
  {\path{doi:{10.1088/0965-0393/18/1/015012}}}.

\bibitem{plimpton1995fast}
S.~Plimpton, Fast parallel algorithms for short-range molecular dynamics,
  Journal of computational physics 117~(1) (1995) 1--19.

\bibitem{wood2018multiscale}
M.~A. Wood, D.~E. Kittell, C.~D. Yarrington, A.~P. Thompson, Multiscale
  modeling of shock wave localization in porous energetic material, Physical
  Review B 97~(1) (2018) 014109.

\bibitem{mishin2001structural}
Y.~Mishin, M.~Mehl, D.~Papaconstantopoulos, A.~Voter, J.~Kress, Structural
  stability and lattice defects in copper: Ab initio, tight-binding, and
  embedded-atom calculations, Physical Review B 63~(22) (2001) 224106.

\bibitem{bringa2004atomistic}
E.~Bringa, J.~Cazamias, P.~Erhart, J.~St{\"o}lken, N.~Tanushev, B.~Wirth,
  R.~Rudd, M.~Caturla, Atomistic shock hugoniot simulation of single-crystal
  copper, Journal of Applied Physics 96~(7) (2004) 3793--3799.

\bibitem{tuckerman2006liouville}
M.~E. Tuckerman, J.~Alejandre, R.~L{\'o}pez-Rend{\'o}n, A.~L. Jochim, G.~J.
  Martyna, A liouville-operator derived measure-preserving integrator for
  molecular dynamics simulations in the isothermal--isobaric ensemble, Journal
  of Physics A: Mathematical and General 39~(19) (2006) 5629.

\bibitem{parrinello1981polymorphic}
M.~Parrinello, A.~Rahman, Polymorphic transitions in single crystals: A new
  molecular dynamics method, Journal of Applied physics 52~(12) (1981)
  7182--7190.

\bibitem{cyr2019robust}
E.~C. Cyr, M.~A. Gulian, R.~G. Patel, M.~Perego, N.~A. Trask, Robust training
  and initialization of deep neural networks: {A}n adaptive basis viewpoint,
  arXiv preprint arXiv:1912.04862 (2019).

\bibitem{you2020data}
H.~You, Y.~Yu, N.~Trask, M.~Gulian, M.~D'Elia, Data-driven learning of robust
  nonlocal physics from high-fidelity synthetic data, arXiv preprint
  arXiv:2005.10076 (2020).

\bibitem{patelaphysics}
R.~G. Patela, N.~A. Traska, M.~A. Woodb, E.~C. Cyra, A physics-informed
  operator regression framework for extracting data-driven continuum models.

\bibitem{Kingma2014}
D.~P. Kingma, J.~Ba, \href{http://arxiv.org/abs/1412.6980}{{Adam: A Method for
  Stochastic Optimization}} (dec 2014).
\newblock \href {http://arxiv.org/abs/1412.6980} {\path{arXiv:1412.6980}}.
\newline\urlprefix\url{http://arxiv.org/abs/1412.6980}

\bibitem{Cyr2020}
E.~C. Cyr, M.~A. Gulian, R.~G. Patel, M.~Perego, N.~A. Trask,
  \href{http://proceedings.mlr.press/v107/cyr20a.html}{Robust training and
  initialization of deep neural networks: {A}n adaptive basis viewpoint}, in:
  J.~Lu, R.~Ward (Eds.), Proceedings of The First Mathematical and Scientific
  Machine Learning Conference, Vol. 107 of Proceedings of Machine Learning
  Research, PMLR, Princeton University, Princeton, NJ, USA, 2020, pp. 512--536.
\newline\urlprefix\url{http://proceedings.mlr.press/v107/cyr20a.html}

\bibitem{Glorot2010}
X.~Glorot, Y.~Bengio,
  \href{http://proceedings.mlr.press/v9/glorot10a.html}{Understanding the
  difficulty of training deep feedforward neural networks}, in: Y.~W. Teh,
  M.~Titterington (Eds.), Proceedings of the Thirteenth International
  Conference on Artificial Intelligence and Statistics, Vol.~9 of Proceedings
  of Machine Learning Research, JMLR Workshop and Conference Proceedings, Chia
  Laguna Resort, Sardinia, Italy, 2010, pp. 249--256.
\newline\urlprefix\url{http://proceedings.mlr.press/v9/glorot10a.html}

\end{thebibliography}

\appendix

\section{Hyper-parameters}

This section lists the hyperparameters used to generate our results. For all neural network architectures for the solutions to PDEs, we use densely connected neural networks of width 64 and depth 8. For neural networks for the EOS, we use densely connected neural networks of width 4 and depth 4 with tanh activation functions. We use the Adam optimizer \cite{Kingma2014} for all minimization problems.

\begin{table}[h]
\begin{center}
\begin{tabular}{|l|l|} 
\hline
Network Initialization & Box \cite{Cyr2020} \\ \hline
TVD hyperparameter ($\epsilon_{TVD}$) & 0 \\ \hline
Activation Function & Elu \\ \hline 
Learning Rate & 1e-4 \\ \hline
Training steps & 4.5e6 \\ \hline
Cells along time dimension & 25 \\ \hline
Cells along space dimension & 50 \\ \hline
Quadrature scheme & Composite trapezoidal \\ \hline
Quadrature segments & 3 \\ 
\hline
\end{tabular}
\end{center}
\caption{Parameters used for the Burgers rarefaction Riemann problem  in Figure \ref{fig:ent}.}\label{table:ent_param}
\end{table}

\begin{table}[h]
\begin{center}
\begin{tabular}{|l|l|} 
\hline
Network Initialization & Box \\ \hline
Entropy hyperparameter ($\epsilon_{ent}$) & 1 \\ \hline
Activation Function & Relu \\ \hline 
Learning Rate & 1e-3 \\ \hline
Training Steps & 9e4 \\ \hline
Cell along time dimension & 64 \\ \hline
Cell along space dimension & 64 \\ \hline
Quadrature scheme & Composite trapezoidal \\ \hline
Quadrature segments & 3 \\ 
\hline
\end{tabular}
\end{center}
\caption{Parameters used for the Sod shock problem in Figure \ref{fig:tvd}.}\label{table:tvd_param}
\end{table}

\begin{table}[h]
\begin{center}
\begin{tabular}{|l|l|} 
\hline
Network Initialization & Glorot \cite{Glorot2010} \\ \hline
Activation Function & ReLU \\ \hline 
Learning Rate & 1e-4 \\ \hline
Training steps & 1e6 \\ \hline
Cells along time dimension & 200 \\ \hline
Cells along space dimension & 200 \\ \hline
Quadrature scheme & Midpoint \\ \hline
Quadrature segments & 1 \\ 
\hline
\end{tabular}
\end{center}
\caption{Parameters used for Buckley-Leverett and Euler Riemann problems in Figure~\ref{fig:forward}.}
\label{tab:gen}
\end{table}

\begin{table}[h]
\begin{center}
\begin{tabular}{|l|l|} 
\hline
Network Initialization & Box \\ \hline
Activation Function & Elu \\ \hline 
Learning Rate & 1e-5 \\ \hline
Training steps & 9e4 \\ \hline
Cell volumes along time dimension & 25 \\ \hline
Cell volumes along space dimension & 50 \\ \hline
Quadrature scheme & Composite trapezoidal \\ \hline
Quadrature segments & 3 \\ 
\hline
\end{tabular}
\end{center}
\caption{Parameters used for the cvPINNs solution of the Burgers shock problem in Figure \ref{fig:PINNs}.}
\label{table:pinn_compare_param}
\end{table}

\begin{table}[h]
\begin{center}
\begin{tabular}{|l|l|} 
\hline
Network Initialization & Glorot \\ \hline
Activation Function & tanh \\ \hline 
Learning Rate & 1e-5 \\ \hline
Training steps & 9e4 \\ \hline
Points along time dimension & 25 \\ \hline
Points along space dimension & 50 \\ \hline
\end{tabular}
\end{center}
\caption{Parameters used for the PINNs solution of the Burgers shock problem in Figure \ref{fig:PINNs}.}
\label{table:pinn_compare_param2}
\end{table}

\begin{table}[h]
\begin{center}
\begin{tabular}{|l|l|} 
\hline
Network Initialization & Glorot \\ \hline
Activation Function & ReLU \\ \hline 
Learning Rate, Phase 1 & 1e-2 \\ \hline
Learning Rate, Phase 2 & 1e-3 \\ \hline
Training steps, Phase 1 & 1e3 \\ \hline
Training steps, Phase 2 & 2e4 \\ \hline
Cells along time dimension & 200 \\ \hline
Cells along space dimension & 200 \\ \hline
Quadrature scheme & Midpoint \\ \hline
Quadrature segments & 1 \\ \hline
Thermodynamic penalty & 100 \\ \hline
Data penalty & 0.1 \\ \hline
\end{tabular}
\end{center}
\caption{Parameters used in Section~\ref{sec:dsmc} and \ref{sec:lammps}. Two training phases with different learning rates and steps were used for these studies.}
\label{tab:inv}
\end{table}

\section{Data sets}

This section lists the parameters used to perform the DSMC simulations of the Sod shock problem and the LAMMPS simulations of the Copper impact problem.

\begin{table}[h]
\begin{center}
\begin{tabular}{|l|l|l|l|l|} 
\hline
Case & $\rho_l$ & $\rho_r$ & $T_l$ & $T_r$ \\
\hline
1 & 1.1 & 0.5 & 4805 & 3844 \\
2 & 0.8 & 0.2 & 1922 & 1922 \\
3 & 1.2 & 0.7 & 7047 & 5766 \\
4 & 0.8 & 0.5 & 5125 & 2562 \\
5 & 1.3 & 0.7 & 3523 & 3523 \\
\hline
\end{tabular}
\end{center}
\caption{Parameters used for the DSMC simulations. The units for density and temperature are kg/m$^3$ and K, respectively.}
\label{tab:dsmc}
\end{table}

\begin{table}[h]
\begin{center}
\begin{tabular}{|l|l|l|} 
\hline
Case & $v_0$ & $T_0$ \\
\hline
1 & 1.0 & 300 \\
2 & 1.5 & 300 \\
3 & 2.0 & 300 \\
4 & 1.0 & 600 \\
5 & 1.5 & 600 \\
6 & 2.0 & 600 \\
7 & 1.0 & 1000 \\
8 & 1.5 & 1000 \\
9 & 2.0 & 1000 \\
\hline
\end{tabular}
\end{center}
\caption{Parameters used for the LAMMPS simulations. The units for velocity and temperature are km/s and K, respectively.}
\label{tab:lammps}
\end{table}

\section{Finite difference scheme} \label{app:fd}

To produce solutions to the Euler equations, Eqn.~\ref{eq:euler}, with the learned EOS's, we use a viscously regularized center difference scheme on a very fine mesh. We define a regular Cartesian grid with points $\{x_i\}$ and grid spacing, $\Delta x$. To solve the Euler equations at these grid points and times, $\{t_n\}$, with timesteps, $\Delta t$, we use the update rule,
\begin{equation}
    \begin{split}
    \bm{u}_{i,n+1} = \bm{u}_{i,n} + \frac{1}{2 \Delta x} (\bm{F}(\bm{u})_{i+1,n}-\bm{F}(\bm{u})_{i-1,n}) \\
    +  \frac{\nu}{\Delta x^2} (\bm{u}_{i,n+1}-2\bm{u}_{i,n} +\bm{u}_{i,n-1})
\end{split}
\end{equation}
where $\nu$ is an artificial viscosity. We use a grid of size 4000, $\nu = 10^{-3}$, $\Delta t=\frac{\Delta x}{32}$. We found this scheme and parameters to produce stable and sharp solutions for the problems considered here.

\end{document}